\documentclass{article}

\usepackage{amssymb}

\usepackage{proof}
\usepackage{latexsym}

\newcommand{\sfZF}{{\sf ZF}}

\newcommand{\KP}{{\sf KP}}

\newcommand{\dg}{\mbox{{\rm dg}}}

\newcommand{\blem}{\begin{lemma}}
\newcommand{\elem}{\end{lemma}}
\newcommand{\bth}{\begin{theorem}}
\newcommand{\ethm}{\end{theorem}}
\newcommand{\benu}{\begin{enumerate}}
\newcommand{\eenu}{\end{enumerate}}
\newcommand{\bdes}{\begin{description}}
\newcommand{\edes}{\end{description}}
\newcommand{\bdf}{\begin{definition}}
\newcommand{\edf}{\end{definition}}
\newcommand{\bcor}{\begin{cor}}
\newcommand{\ecor}{\end{cor}}
\newcommand{\bprp}{\begin{proposition}}
\newcommand{\eprp}{\end{proposition}}
\newcommand{\bmlem}{\begin{mlemma}}
\newcommand{\emlem}{\end{mlemma}}
\newcommand{\bclm}{\begin{claim}}
\newcommand{\eclm}{\end{claim}}
\newcommand{\brem}{\begin{remark}}
\newcommand{\erem}{\end{remark}}
\newcommand{\bprf}{{\bf Proof}.\hspace{2mm}}
\newcommand{\eprf}{\hspace*{\fill} $\Box$}

\newcommand{\beqn}{\begin{equation}}
\newcommand{\eeqn}{\end{equation}}
\newcommand{\beqnarr}{\begin{eqnarray}}
\newcommand{\eeqnarr}{\end{eqnarray}}
\newcommand{\beqnarrs}{\begin{eqnarray*}}
\newcommand{\eeqnarrs}{\end{eqnarray*}}
\newcommand{\spand}{\,\&\,}

\newcommand{\restrict}{\!\upharpoonright\!}

\newtheorem{theorem}{Theorem}[section]
\newtheorem{definition}[theorem]{Definition}
\newtheorem{proposition}[theorem]{Proposition}
\newtheorem{lemma}[theorem]{Lemma}
\newtheorem{cor}[theorem]{Corollary}
\newtheorem{claim}[theorem]{Claim}
\newtheorem{remark}[theorem]{Remark}

\newcommand{\alp}{\alpha}

\newcommand{\veps}{\varepsilon}
\newcommand{\del}{\delta}
\newcommand{\Del}{\Delta}
\newcommand{\ome}{\omega}
\newcommand{\Ome}{\Omega}
\newcommand{\bet}{\beta}
\newcommand{\gam}{\gamma}
\newcommand{\Gam}{\Gamma}
\newcommand{\kap}{\kappa}
\newcommand{\sig}{\sigma}
\newcommand{\Sig}{\Sigma}
\newcommand{\tht}{\theta}

\newcommand{\lam}{\lambda}
\newcommand{\Lam}{\Lambda}
\newcommand{\vphi}{\varphi}

\newcommand{\fal}{\forall}
\newcommand{\exi}{\exists}

\newcommand{\Rarw }{\Rightarrow}

\newcommand{\lrarw}{\leftrightarrow}
\newcommand{\Lrarw}{\Leftrightarrow}

\newcommand{\cala}{{\cal A}}

\newcommand{\calh}{{\cal H}}

\newcommand{\calL}{{\cal L}}

\newcommand{\calP}{{\cal P}}

\newcommand{\la}{\langle}
\newcommand{\ra}{\rangle}

\newcommand{\msten}{\mbox{\hspace{10mm}}}
\newcommand{\msfiv}{\mbox{\hspace{5mm}}}

\title{
Cut-elimination for $\ome_{1}$
}

\author{Toshiyasu Arai
\\
Graduate School of Science,
Chiba University
\\
1-33, Yayoi-cho, Inage-ku,
Chiba, 263-8522, JAPAN
\\
tosarai@faculty.chiba-u.jp
}

\begin{document}
\maketitle

\begin{abstract}
In this paper 
we calibrate the strength of the soundness of a set theory $\KP\ome+(\Pi_{1}\mbox{{\rm -Collection}})$
with the assumption that
`there exists an uncountable regular ordinal'
in terms of the existence of ordinals.
\end{abstract}

\section{Introduction}\label{sect:intro}

In \cite{liftupK, liftupZF} higher set theories are analyzed proof-theoretically
in terms of the operator controlled derivations, which are introduced by W. Buchholz\cite{Buchholz92}.
Let $I$ be the least weakly inaccessible cardinal.
Collapsing functions $\alp\mapsto\Psi_{\kap,n}(\alp)<\kap$ are introduced 
for each uncountable regular cardinal $\kap\leq I$ and $n<\ome$.
$\Psi_{\kap,n}(\alp)$ is a first-order variant of the collapsing functions $\psi_{\kap}(\alp)$ 
introduced by W. Buchholz\cite{Buchholz86}.
Let $\ome_{k}(I+1)$ denote the tower of $\ome$ with the next epsilon number 
$\veps_{I+1}=\sup\{\ome_{k}(I+1): k<\ome\}$ above $I$.
The predicate $x=\Psi_{\kap,n}(\alp)$ is a $\Sig_{n+1}$-predicate for $\alp<\veps_{I+1}$, and
we see that for each $n,k<\ome$
$\sfZF+(V=L)$ proves $\fal\alp<\ome_{k}(I+1)
\fal\kap\leq I\exi x<\kap[x=\Psi_{\kap,n}(\alp)]$.

Conversely the following Theorem \ref{th:mainthZ} is shown in \cite{liftupZF}.

\bth\label{th:mainthZ}
For a sentence $\exi x\in L_{\ome_{1}}\, \vphi(x)$ with a first-order formula $\vphi(x)$, 
if
\[
\sfZF+(V=L)\vdash\exi x\in L_{\ome_{1}}\,\vphi(x) 
\]
then
\[
\exi n<\ome[\sfZF+(V=L)\vdash \exi x\in L_{\Psi_{\ome_{1},n}(\ome_{n}(I+1))} \vphi(x)]
.\]
\end{theorem}
Theorem \ref{th:mainthZ} is shown as follows.
First the finite $\sfZF$-proofs are embedded into infinitary operator controlled derivations.
Second cut inferences are eliminated from the infinitary derivations.
Third we conclude that the end formula is true by transfinite induction on the depths of the cut-free derivations.
To formalize the proof we need a derivability predicate
$\calh_{\gam}\vdash^{\alp}_{c}\Gam$ for operator controlled derivations,
the existence of the collapsing functions $\Psi_{\kap,n}(\alp)$ and
the transfinite induction $TI(\veps_{I+1})$.
We can define the predicate $\calh_{\gam}\vdash^{\alp}_{c}\Gam$ as a fixed point,
and the argument with respect to the predicate is carried out intuitionistically.

For \textit{each} $\sfZF$-proof, we can find a $k<\ome$ such that the transfinite induction
$TI(\ome_{k}(I+1))$ suffices to formalize the proof.
Hence the whole proof is formalized in
an intuitionistic fixed point theory ${\rm Fix}^{i}(\sfZF)$ over $\sfZF$,
which is a conservative extension of $\sfZF$, cf.\,\cite{IntSet}.
Theorem \ref{th:mainthZ} aims at enlarging the realm of the ordinal analysis, a topic in proof theory.

On the other side, theorems in ordinal analysis are typically stated as follows.
Let $T$ be a theory including the first-order arithmetic {\sf PA}.
A computable notation system $O(T)$ of ordinals is defined for the proof-theoretic ordinal of the theory $T$,
for which the following holds.
For example when $T$ is the second-order arithmetic ${\rm ATR}_{0}$ for the arithmetical transfinite recursion,
$O({\rm ATR}_{0})=\Gam_{0}$, the first strongly critical ordinal, and
$O(\Pi^{1}_{1}{\rm -CA}_{0})=\psi_{\Ome_{1}}(\Ome_{\ome})$,
where $\Ome_{\ome}=\sup\{\Ome_{1+n}:n<\ome\}$
and $\Ome_{n}$ is the $n$-th recursively regular ordinal with $\Ome_{1}=\ome_{1}^{CK}$.
$TI(O(T))$ [$TI(<\!\!O(T))$]
denotes the $\Pi^{1}_{1}$-sentence expressing the wellfoundedness of $O(T)$
[the $\Pi^{1}_{1}$-sentences expressing the wellfoundedness up to each ordinal in $O(T)$], resp.

\bth\label{th:intro}
\benu
\item
The uniform reflection principle $\mbox{{\rm RFN}}_{\Pi^{1}_{1}}(T)$, i.e., 
$\Pi^{1}_{1}$-soundness of $T$
is provable from $\mbox{{\rm TI}}(O(T))$ over  {\sf PA}.

\item 
Each $\Pi^{1}_{1}$-provable sentence in $T$ follows from $\mbox{{\rm TI}}(<\! O(T))$ over {\sf PA}.
\item 
$T$ proves $\mbox{{\rm TI}}(<\!O(T))$.
\eenu
\end{theorem}

From this we conclude the following.

\bcor\label{cor:introa}
\benu
\item
$\mbox{{\rm RFN}}_{\Pi^{1}_{1}}(T)$ is equivalent to $\mbox{{\rm TI}}(O(T))$ over a weak theory, e.g., 
over the elementary arithmetic {\sf EA}.

\item
Two theories $T$ and $\mbox{{\rm TI}}(<\!O(T))$ have the same provable $\Pi^{1}_{1}$-theorems.
\eenu
\ecor
Thus the soundness of $T$ (with respect to a class of formulas) is equivalent to
the existence (wellfoundedness) of the ordinal $O(T)$,
and the consequences in $T$ follow from the existence of ordinals$<O(T)$.

In this paper we show a similar result for a higher set theory.
We calibrate the strength of the soundness of a set theory $T_{1}$
in terms of the existence of ordinals, where
$T_{1}:=\KP\ome+(V=L)+(\Pi_{1}\mbox{{\rm -Collection}})+(\ome_{1})$,
and $(\ome_{1})$ denotes an axiom stating that `there exists an uncountable regular ordinal'.
Let $\rho_{0}>\ome_{1}$ be the least ordinal such that $L_{\rho_{0}}\models(\Pi_{1}\mbox{{\rm -Collection}})$.
Let $<^{\veps}$ be a $\Del_{1}$-well ordering whose order type is the next epsilon number
$\veps_{\rho_{0}+1}$ to the order type $\rho_{0}$ of the class of ordinals in a transitive and wellfounded
model of $T_{1}$.
$\alp<\veps_{\rho_{0}+1}$ denotes the fact that $\alp$ is in the field of the ordering $<^{\veps}$.
As in \cite{liftupZF} collapsing functions $\Psi_{\kap}(\alp)<\kap\,(\kap\in\{\ome_{1},\rho_{0}\})$
and the Skolem hulls $\calh_{\alp}(\bet)\,(\bet<\rho_{0})$ are introduced.
Each of $x=\Psi_{\kap}(\alp)$ and $x=\calh_{\alp}(\bet)$ is a $\Sig_{2}$-predicate.

For $\alp\in\{\veps_{\rho_{0}+1}\}\cup\{\ome_{k}(\rho_{0}+1): k<\ome\}$,
$T(\alp)$ denotes a set of ordinal terms $a$ representing ordinals $v(a)$ in $\calh_{\alp}(0)$.
$T(\ome_{k}(\rho_{0}+1))$ is a subset of $T(\veps_{\rho_{0}+1})$.
$T(\veps_{\rho_{0}+1})$ is a computable set of integers under a suitable encoding.
We see that the relation $v(a)=\alp$ for $a\in T(\veps_{\rho_{0}+1})$ and $\alp<\veps_{\rho_{0}+1}$
is a $\Sig_{2}$-predicate in $\KP\ome$, cf.\,Proposition \ref{prp:value}.
Let for $a\in T(\veps_{\rho_{0}+1})$
\beqnarr
&&
A(a)
:\Lrarw
\exi\alp<\veps_{\rho_{0}+1}\{
v(a)=\alp \land
\nonumber
\\
&&
\exi\bet<\rho_{0}[\Psi_{\rho_{0}}(\alp)=\bet] \land 
\exi\bet<\ome_{1}[\Psi_{\ome_{1}}(\alp)=\bet\land F_{\bet\cup\{\ome_{1}\}}(\rho_{0})<\ome_{1}]
\}
\label{eq:A}
\eeqnarr
where
$F_{x\cup\{\ome_{1}\}}(y)$ denotes the Mostowski collapse
 $F_{x\cup\{\ome_{1}\}}:\mbox{{\rm Hull}}(x\cup\{\ome_{1}\})\lrarw \gam$
 of the $\Sig_{1}$-Skolem hull $\mbox{{\rm Hull}}(x\cup\{\ome_{1}\})$ of $x\cup\{\ome_{1}\}$
 in $\rho_{0}$ with
$x=F_{x\cup\{\ome_{1}\}}(\ome_{1})$ and $\gam=F_{x\cup\{\ome_{1}\}}(\rho_{0})$.
$F_{\bet\cup\{\ome_{1}\}}(\rho_{0})=\gam$ denotes a $\Sig_{2}$-predicate such that
$\exi h[h=\mbox{{\rm Hull}}(\bet\cup\{\ome_{1}\})\land \fal\del \in h(F_{\bet\cup\{\ome_{1}\}}(\del)<\gam) \land
\fal\bet<\gam\exi \del\in h(\bet\leq F_{\bet\cup\{\ome_{1}\}}(\del))]$.

$\mbox{{\rm RFN}}_{\Sig_{2}}(T_{1})$ denotes the uniform reflection principle for $T_{1}$ with respect to
$\Sig_{2}$-formulas $\vphi(n)$ with the individual constant $\ome_{1}$ and a variable $n$ for integers:
\[
\fal n<\ome\left[
{\rm Pr}_{T_{1}}(\lceil\vphi(\dot{n})\rceil) \to \vphi(n)
\right]
\]
with a standard provability predicate ${\rm Pr}_{T_{1}}$ for 
$T_{1}$.

$\KP\ome+(V=L)$ denotes the set theory obtained from
 the Kripke-Platek set theory $\KP\ome$ with the axiom of infinity by adding the axiom $V=L$ of constructibility
and the axiom $\ome<\ome_{1}(<\rho_{0})$ for the constant $\ome_{1}$.

\bth\label{th:main}
\benu
\item\label{th:main1}
The uniform reflection principle $\mbox{{\rm RFN}}_{\Sig_{2}}(T_{1})$ is equivalent to
$\fal a\in T(\veps_{\rho_{0}+1})\, A(a)$ over $\KP\ome+(V=L)$.

\item\label{th:main2}
For any $\Sig_{2}$-formula $\vphi(n)$ with parameters $n<\ome$,
$T_{1}\vdash\fal n<\ome\,\vphi(n)$ iff
there exists a $k<\ome$ for which
$\KP\ome+(V=L)\vdash \fal a\in T(\ome_{k}(\rho_{0}+1))\, A(a)\to \fal n<\ome\,\vphi(n)$ holds.
\eenu

\end{theorem}

Similarly for Theorem \ref{th:mainthZ},
 we could show Theorem \ref{th:main}.\ref{th:main2} through an intuitionistic fixed point theory
and the controlled derivations.
However if we follow this tactics in showing a soundness of the set theory $T_{1}$,
we obtain only a weaker assertion:
the uniform reflection principle $\mbox{{\rm RFN}}_{\Sig_{2}}(T_{1})$ follows from
$\fal a\in T(\veps_{\rho_{0}+1})\, A(a)$ and the \textit{transfinite induction} $TI(\veps_{\rho_{0}+1})$.
Apparently $\mbox{{\rm RFN}}_{\Sig_{2}}(T_{1})$ does not yield $TI(\veps_{\rho_{0}+1})$.
Hence we need an alternative treatment to show Theorem \ref{th:main}.\ref{th:main1}.

Our proof is in the scheme
of the consistency proofs in G. Gentzen\cite{Gentzen38} and in G. Takeuti\cite{Takeuti67},
in which ordinals $o(\calP)$ are associated with proof figures $\calP$ in such a way that
$o(\calP)>o(r(\calP))$ for a proof figure $r(\calP)$ of a contradiction.
Given a proof figure $\calP_{0}$ of a contradiction,
define proof figures $\calP_{n}$ of a contradiction recursively by $\calP_{n+1}=r(\calP_{n})$. 
Then the series $\{\alp_{n}\}_{n}$ of ordinals $\alp_{n}=o(\calP_{n})$ would be an infinite descending chain, and
hence there is no proof figures of a contradiction.
We assign ordinals $o(\calP)$ to finite proof figures $\calP$ of $\Sig_{2}$-sentences,
and define a rewriting step $r(\calP)$ on such proof figures $\calP$
 in which constants for ordinals$<\rho_{0}$ may occur.
In \cite{Gentzen38, Takeuti67} both of rewriting step $r$ and ordinal assignment $o$
are primitive (or even elementary) recursive.
In our proof a transcendence over finite mathematics gets into the definition of rewriting steps 
(and the satisfaction relation for $\Del_{0}$-formulas).
\\

\noindent
Let us mention the contents of the paper.
In section \ref{sect:Clps}
let us recall $\Sig_{1}$-Skolem hulls, a paraphrase of the regularity of ordinals, and 
ordinals for regular ordinals.
All of these come from  \cite{liftupZF}.
In section \ref{sect:finiteproof}, 
an ordinal assignment $o(\Gam)$ to sequents $\Gam$ occurring in proofs are defined.
Finally we define a rewriting step $\mathcal{P}\mapsto\mathcal{P}^{\prime}$ on (finite) proof figures for which
$o(\calP)>o(r(\calP))$ holds, and a proof of Theorem \ref{th:main} is concluded in section \ref{sect:consisprf}.

\section{Collapsing functions for $\ome_{1}$}\label{sect:Clps}
In this section
let us recall $\Sig_{1}$-Skolem hulls, a paraphrase of the regularity of ordinals, and 
ordinals for regular ordinals.
Everything in this section is reproduced from \cite{liftupZF}.

\subsection{$\Sig_{1}$-Skolem hulls}\label{sect:Skolemhull}

Let $L_{\alp}$ be the $\alp$-th level of the conctructible universe $L$.
$\rho_{0}$ denotes the least ordinal above $\ome_{1}$ such that 
$L_{\rho_{0}}\models (\Pi_{1}\mbox{{\rm -Collection}})$.

\bdf\label{df:crdS} 

\benu

\item 
$
cf(\kap) := \min\{\alp\leq\kap : \mbox{{\rm there is a cofinal map }} f:\alp\to \kap\}
$.

\item
{\rm For} $X\subset L_{\rho_{0}}$,
$\mbox{{\rm Hull}}(X)$ {\rm denotes the set} 
{\rm (}$\Sig_{1}$-Skolem hull {\rm of} $X$ {\rm in} $L_{\rho_{0}}${\rm ) defined as follows.}
$<_{L}$ {\rm denotes a} $\Del_{1}${\rm -well ordering of the constructible universe} $L$.
{\rm Let} $\{\vphi_{i}:i\in\ome\}$ {\rm denote an enumeration of} $\Sig_{1}${\rm -formulas in the language}
$\{\in\}${\rm . Each is of the form} $\vphi_{i}\equiv\exi y\theta_{i}(x,y;u)\, (\tht\in\Del_{0})$ {\rm with fixed variables} $x,y,u${\rm . Set for} $b\in X$
\beqnarrs
r_{\Sig_{1}}^{\rho_{0}}(i,b) & \simeq & \mbox{ {\rm the }} <_{L} \mbox{{\rm -least }} c\in L_{\rho_{0}}
\mbox{ {\rm such that} } L_{\rho_{0}}\models\theta_{i}((c)_{0},(c)_{1}; b)
\\
h_{\Sig_{1}}^{\rho_{0}}(i,b) & \simeq & (r_{\Sig_{1}}^{\rho_{0}}(i,b))_{0}
\\
\mbox{{\rm Hull}}(X) & = & rng(h_{\Sig_{1}}^{\rho_{0}})=\{h_{\Sig_{1}}^{\rho_{0}}(i,b)\in L_{\rho_{0}}:i\in\ome, b\in X\}
\eeqnarrs
{\rm Then}
$
L_{\rho_{0}}\models \exi x\exi y\, \theta_{i}(x,y;b) \to h_{\Sig_{1}}^{\rho_{0}}(i,b)\downarrow \spand \exi y\, \theta_{i}(h_{\Sig_{1}}^{\rho_{0}}(i,b),y;b)
$.

\item
{\rm The Mostowski collapsing function}
\[
F_{X}:\mbox{{\rm Hull}}(X)\lrarw  L_{\gam}
\]
{\rm for an ordinal} $\gam\leq\rho_{0}$ {\rm such that} $F_{X}\restrict Y=id\restrict Y$ 
{\rm for any transitive}
$Y\subset \mbox{{\rm Hull}}(X)$.

{\rm Let us denote, though} $\rho_{0}\not\in dom(F)=\mbox{{\rm Hull}}(X)$
\[
F_{X}(\rho_{0}):=\gam
.\]

\eenu
\edf

\bprp\label{prp:complexcr}
Let 
$L_{\rho_{0}}\models\KP\ome$.
Then for $\kap\leq\rho_{0}$,
$\{(x,y): x<\kap\spand y=\min\{y<\kap: \mbox{{\rm Hull}}(x\cup\{\kap\})\cap\kap\subset y\}\}$ 
is a $Bool(\Sig_{1}(L_{\rho_{0}}))$-predicate on $\kap$, and hence the set is in $L_{\rho_{0}}$
if $\kap<\rho_{0}$ and $L_{\rho_{0}}\models\Sig_{1}\mbox{{\rm -Separation}}$.
\eprp

\bprp\label{prp:rhHull}
Assume that $X$ is a set in $L_{\rho_{0}}$.
Then
$r_{\Sig_{1}}^{\rho_{0}}$ and $h_{\Sig_{1}}^{\rho_{0}}$ are partial $\Del_{1}(L_{\rho_{0}})$-maps such that the domain of $h_{\Sig_{1}}^{\rho_{0}}$ is a 
$\Sig_{1}(L_{\rho_{0}})$-subset of 
$\ome\times X$.
Therefore its range $\mbox{{\rm Hull}}(X)$ is a $\Sig_{1}(L_{\rho_{0}})$-subset of $L_{\rho_{0}}$.
\eprp

\bth\label{th:cofinalitylocal}
Let $\rho_{0}$ be an ordinal such that $L_{\rho_{0}}\models \mbox{{\rm KP}}\ome+(\Pi_{1}\mbox{{\rm -Collection}})$,
 and $\ome\leq\alp<\kap<\rho_{0}$ with
$\alp$ a multiplicative principal number and $\kap$ a limit ordinal.
Then the following conditions are mutually equivalent:
\benu

\item 
$L_{\rho_{0}}\models {}^{\alp}\kap\subset L_{\kap}$.

\item 
$L_{\rho_{0}}\models \alp<cf(\kap)$.

\item
There exists an ordinal $x$ such that $\alp<x<\kap$,
$ \mbox{{\rm Hull}}(x\cup\{\kap\})\cap \kap\subset x$ and $F_{x\cup\{\kap\}}(\rho_{0})<\kap$.

\item 
For the Mostowski collapse
$F_{x\cup\{\kap\}}(y)$, there exists an ordinal $x$ such that 
$\alp<x=F_{x\cup\{\kap\}}(\kap)<F_{x\cup\{\kap\}}(\rho_{0})<\kap$,
and for any $\Sig_{1}$-formula $\vphi$ and any $a\in L_{x}$,
$L_{\rho_{0}}\models\vphi[\kap,a]\Rarw  L_{F_{x\cup\{\kap\}}(\rho_{0})} \models\vphi[x,a]$ holds.

\eenu

\end{theorem}

\subsection{Theories equivalent to $T_{1}$}\label{sect:Ztheory}

Referring Theorem \ref{th:cofinalitylocal} we introduce a theory $T(\ome_{1})$ equivalent to $T_{1}$.

\bdf\label{df:regext}
$T(\ome_{1})$ {\rm denotes the set theory defined as follows.}
\benu
\item
{\rm Its language is} $\{\in, P,P_{\rho_{0}},\ome_{1}\}$ {\rm for a binary predicate} $P$, {\rm a unary predicate} $P_{\rho_{0}}$
 {\rm  and an individual constant} $\ome_{1}$.

\item
{\rm Its axioms are obtained from those of} $\KP\ome+(\Pi_{1}\mbox{{\rm -Collection}})$ 
{\rm
 in the expanded language\footnote{
This means that the predicates $P,P_{\rho_{0}}$ do not occur in 
$\Del_{0}$-formulas
for $\Del_{0}$-Separation and 
$\Pi_{1}$-formulas
for $\Pi_{1}$-Collection.
},
the axiom of constructibility}
$V=L$
{\rm together with the axiom schema saying that}
$\ome_{1}$ 
 {\rm is an uncountable regular ordinal, cf.\,(\ref{eq:Z2}) and (\ref{eq:Z1}),
 and if} $P(x,y)$ {\rm then}  $x$ 
 {\rm is a critical point of the} $\Sig_{1}${\rm -elementary embedding from}
$L_{y}\cong \mbox{{\rm Hull}}(x\cup\{\ome_{1}\})$ {\rm to the universe}
$L_{\rho_{0}}${\rm , cf.\,(\ref{eq:Z1}), and if}
$P_{\rho_{0}}(x)$ {\rm then} $x$ {\rm is a critical point of the} $\Sig_{1}${\rm -elementary embedding from}
$L_{x}\cong \mbox{{\rm Hull}}(x)$ {\rm to the universe}
$L_{\rho_{0}}${\rm , cf.\,(\ref{eq:Z4}):}
{\rm for a formula} $\vphi$ {\rm and an ordinal} $\alp$,
$\vphi^{\alp}$ {\rm denotes the result of restricting every unbounded quantifier}
$\exi z,\fal z$ {\rm in} $\vphi$ {\rm to} $\exi z\in L_{\alp}, \fal z\in L_{\alp}$.

 \benu
 
 \item
 $x\in Ord$ {\rm is a} $\Del_{0}${\rm -formula saying that `}$x$ {\rm is an ordinal'.}
 \\
$ (\ome<\ome_{1}\in Ord)$, $(P(x,y) \to \{x,y\}\subset Ord \land  x<y<\ome_{1})$
{\rm and}
$(P_{\rho_{0}}(x) \to x\in Ord)$.

 \item
\beqn\label{eq:Z1}
P(x,y) \to a\in L_{x}  \to \vphi[\ome_{1},a] \to \vphi^{y}[x,a]
\eeqn
{\rm for any} $\Sig_{1}${\rm -formula} $\vphi$ {\rm in the language} $\{\in\}$.

.

\item
\beqn\label{eq:Z2}
a\in Ord\cap\ome_{1} \to \exi x, y\in Ord\cap\ome_{1}[a<x\land P(x,y)]
\eeqn

\item
\beqn\label{eq:Z4}
P_{\rho_{0}}(x) \to a\in L_{x} \to \vphi[a] \to \vphi^{x}[a]
\eeqn
{\rm for any} $\Sig_{1}${\rm -formula} $\vphi$ {\rm in the language} $\{\in\}$.

\item
\beqn\label{eq:Z5}
a\in Ord \to \exi x\in Ord[a<x\land P_{\rho_{0}}(x)]
\eeqn

 \eenu
 
\eenu
\edf

\brem
{\rm
Though the axioms (\ref{eq:Z4}) and (\ref{eq:Z5}) for the $\Pi_{1}$-definable predicate $P_{\rho_{0}}(x)$ are derivable
from $\Pi_{1}$-Collection, the primitive predicate symbol $P_{\rho_{0}}(x)$ is useful for our proof-theoretic study, 
cf. {\bf Case 1} in subsection \ref{subsec:toprule}.
}
\erem

\blem\label{lem:regularset}
$T(\ome_{1})$ is 
a 
conservative extension of
the set theory $T_{1}$.
\elem
\bprf
First consider the axioms of $T_{1}$ in $T(\ome_{1})$.
The axiom $(\ome_{1})$ is codified as
\[
(\ome_{1}) \:\:  \exi\kap\in Ord[\ome<\kap=cf(\kap)]
\]
which follows from (\ref{eq:Z1}) in $T(\ome_{1})$.
Hence $T_{1}$ is a subtheory of $T(\ome_{1})$.

Next we show that $T(\ome_{1})$ is interpretable in $T_{1}$.
Let $\kap$ be an ordinal in the axiom $(\ome_{1})$.
Interpret the predicate
$P(x,y)\lrarw  (\{x,y\}\subset Ord) \land 
(\mbox{{\rm Hull}}(x\cup\{\kap\})\cap \kap\subset x)
\land
(y=\sup\{F_{x\cup\{\kap\}}(a):a\in\mbox{{\rm Hull}}(x\cup\{\kap\})\})$.
We see from Theorem \ref{th:cofinalitylocal} that
the interpreted (\ref{eq:Z1}) and (\ref{eq:Z2})
are provable in $T_{1}$.

It remains to show the interpreted (\ref{eq:Z4}) and (\ref{eq:Z5}) in $T_{1}$.
It suffices to show that given an ordinal $\alp$,
there exists an ordinal $x>\alp$
such that $\mbox{Hull}(x)\cap Ord\subset x$.

First we show that for any $\alp$ there exists a $\bet$
such that $\mbox{Hull}(\alp)\cap Ord\subset \bet$.
By Proposition \ref{prp:rhHull}
let $h^{\rho_{0}}_{\Sig_{1}}$ be the $\Del_{1}$-surjection from the $\Sig_{1}$-subset
$dom(h^{\rho_{0}}_{\Sig_{1}})$ of $\ome\times\alp$ to
$\mbox{Hull}(\alp)$, which is 
a $\Sig_{1}$-class.
From $\Sig_{1}$-Separation 
we see that $dom(h^{\rho_{0}}_{\Sig_{1}})$ is a set.
Hence by $\Sig_{1}$-Collection,
$\mbox{Hull}(\alp)=rng(h^{\rho_{0}}_{\Sig_{1}})$ is a set.
Therefore the ordinal $\sup(\mbox{Hull}(\alp)\cap Ord)$ exists
in the universe.

As in Proposition \ref{prp:complexcr}
we see that
$X=\{(\alp,\bet): \bet=\min\{\bet\in Ord: \mbox{{\rm Hull}}(\alp)\cap Ord\subset \bet\}\}$ 
is  a set in $L_{\rho_{0}}$ as follows.
Let $\vphi(\bet)$ be the $\Pi_{1}$-predicate $\vphi(\bet):\Lrarw 
\fal \gam\in Ord[\gam\in\mbox{{\rm Hull}}(\alp)\to \gam\in \bet]$.
Then
$\bet=\min\{\bet: \mbox{{\rm Hull}}(\alp)\cap Ord\subset \bet\}$ iff
$\vphi(\bet)
\land \fal \gam<\bet\lnot \vphi(\gam)$, which is  
$Bool(\Sig_{1}(L_{\rho_{0}}))$ by $\Pi_{0}\mbox{{\rm -Collection}}$.
Hence $X$ is a set in $L_{\rho_{0}}$.

Define recursively ordinals $\{x_{n}\}_{n}$ as follows.
$x_{0}=\alp+1$, and $x_{n+1}$ is defined to be the least ordinal $x_{n+1}$ 
such that
$\mbox{{\rm Hull}}(x_{n})\cap Ord\subset x_{n+1}$, i.e., $(x_{n},x_{n+1})\in X$.
We see inductively that such an ordinal $x_{n}$ exists.
Moreover $n\mapsto x_{n}$ is a $\Del_{1}$-map.
Then $x=\sup_{n}x_{n}<\rho_{0}$ is a desired one.
\eprf
\\

Next let us interpret the set theory $T(\ome_{1})$ in a theory $T^{ord}(\ome_{1})$ of ordinals
as in \cite{ptMahlo}.
The base language is $\mathcal{L}_{0}=\{<,0,+,\cdot,\lam x.\ome^{x}\}$.
Each of functions $1,\max$ and the G\"odel pairing function $j$ is $\Del_{0}$-definable in $\mathcal{L}_{0}$, 
cf.\,Appendix B
of \cite{ptMahlo}.
For each bounded formula $\cala(X,a,b)$ in the base language $\mathcal{L}_{0}$, introduce a binary predicate symbol $R^{\cala}$
with its defining axiom $b\in R^{\cala}_{a}:\equiv R^{\cala}(a,b) \lrarw \cala(R^{\cala}_{<a},a,b)$
where $c\in R^{\cala}_{<a}:\Lrarw \exi d<a(c\in R^{\cala}_{d})$.
$\mathcal{L}_{1}$ denotes the resulting language with these predicates $R^{\cala}$.
$\KP\ome+(V=L)$ is interpretable in a theory $T_{2}$ with the axiom for $\Pi_{2}$-reflection, cf.\,Appendix A
of \cite{ptMahlo}.
Each epsilon number $\alp$ is identified with the $\mathcal{L}_{1}$-structure
$\la\alp; <,0,+,\cdot,\lam x.\ome^{x}, R^{\cala}\ra$.
A G\"odel's surjective map $F:Ord\to L$ maps each epsilon number (or even a multiplicative principal number) $\alp$ onto $L_{\alp}$,
and $a\epsilon b\Lrarw F(a)\in F(b)\,(a,b\in L_{\alp})$ is a $\Del_{0}$-relation in the language $\mathcal{L}_{1}$.

For $\Pi_{2}$-formula $A$ in the language $\mathcal{L}_{1}$,
$A(t) \to \exi y[t<y\land A^{(y)}(t)]$
is an instance of $\Pi_{2}$-reflection, which follows from $(V=L)$ and $\Del_{0}$-Collection,
where $A^{(y)}$ denotes the result of restricting unbounded quantifiers $Q x\,(Q\in\{\exi,\fal\})$ to $Qx<y$.

The language of the theory $T^{ord}(\ome_{1})$ is defined to be 
$\mathcal{L}_{2}=\mathcal{L}_{1}\cup\{\ome_{1},P,P_{\rho_{0}}\}$.
The axiom (\ref{eq:Z1}) is translated to
\beqn\label{eq:Z1ord}
P(x,y) \to a<x \to \vphi[\ome_{1},a] \to \vphi^{y}[x,a]
\eeqn
 for $\Sig_{1}$-formulas $\vphi$ in $\mathcal{L}_{1}$.
The axiom (\ref{eq:Z2}) becomes
\beqn\label{eq:Z2ord}
a<\ome_{1} \to \exi x, y<\ome_{1}[a<x\land P(x,y)]
\eeqn
The axiom (\ref{eq:Z4}) turns to
\beqn\label{eq:Z4ord}
P_{\rho_{0}}(x) \to a<x \to \vphi[a] \to \vphi^{x}[a]
\eeqn
The axiom (\ref{eq:Z5}) is formulated in
\beqn\label{eq:Z5ord}
\exi x[a<x\land P_{\rho_{0}}(x)]
\eeqn

Finally consider $\Pi_{1}$-Collection.
For a $\Del_{0}$-formula $\tht(u,v,w)$ in the language $\{\in\}$,
let $\tau(x,u,v)\equiv[u,v\in L_{x}\land \fal w\in L_{x}\tht(u,v,w)]$.
Then we see $\fal w\,\tht(u,v,w)\lrarw\exi x\in P_{\rho_{0}}\tau(x,u,v)$ from $(V=L)$, (\ref{eq:Z4}) and (\ref{eq:Z5}).
Hence $\Pi_{1}$-Collection
\[
\fal u\in a\exi v \fal w\,\tht \to\exi c[a\in c\land \fal u\in a\exi v\in c\fal w\,\tht]\:
(c \mbox{ is transitive})
\]
follows from
\[
\fal u\in a\, A(u) \to \exi c[a\in c\land \fal u\in a\, A^{(c)}(u)]
\]
where $A(u)\equiv(\exi x\in P_{\rho_{0}}\exi v\, \tau(x,u,v))$ for  
$A^{(c)}(u)\equiv(\exi x\in P_{\rho_{0}}\cap c\exi v\in c\,\tau)$.
The latter is translated in the language $\mathcal{L}_{2}$ to
\beqn\label{eq:Z6ord}
\fal u<a\, A(u) \to \exi c>a\fal u<a\, A^{(c)}(u)
\eeqn
where $A(u)\equiv(\exi x\in P_{\rho_{0}}\exi v\, \tau(x,u,v))$ with a $\Del_{0}$-formula $\tau$ in $\mathcal{L}_{1}$.

Let $\mbox{{\rm T}}^{ord}(\ome_{1})$ denote the resulting extension of the theory $T_{2}$ of ordinals 
with axioms (\ref{eq:Z1ord}), (\ref{eq:Z2ord}), (\ref{eq:Z4ord}), (\ref{eq:Z5ord}) and (\ref{eq:Z6ord}),
in which $\mbox{{\rm T}}(\ome_{1})$ is interpreted.

\subsection{Ordinals for $\ome_{1}$}\label{sect:ordinals}

Let $Ord^{\varepsilon}$ and $<^{\varepsilon}$ be $\Delta$-predicates on the universe $V$ such that
for any transitive and wellfounded model $V$ of $\mbox{{\sf KP}}\omega$,
$<^{\varepsilon}$ is a well ordering of type $\varepsilon_{\rho_{0}+1}$ on $Ord^{\varepsilon}$
for the order type $\rho_{0}$ of the class $Ord$ in $V$.
$<^{\veps}$ is seen to be a canonical ordering as stated in
the following Proposition \ref{prp:canonical}.
For natural numbers $n$,
$\ome_{n}(\rho_{0}+1)\in Ord^{\veps}$ is defined recursively by
$\ome_{0}(\rho_{0}+1)=\rho_{0}+1$ and $\ome_{n+1}(\rho_{0}+1)=\ome^{\ome_{n}(\rho_{0}+1)}$.

\bprp\label{prp:canonical}
\benu
\item\label{prp:canonical0}
$\KP\ome$ proves the fact that $<^{\veps}$ is a linear ordering.

\item\label{prp:canonical2}
For any formula $\vphi$
and each $n<\ome$, $\KP\ome$ proves the transfinite induction schema up to $\ome_{n}(\rho_{0}+1)$,
$
\fal x\in Ord^{\veps}(\fal y<^{\veps}x\,\vphi(y)\to\vphi(x)) \to 
\fal x<^{\veps}\ome_{n}(\rho_{0}+1)\vphi(x)$.
\eenu
\eprp


For simplicity let us identify the code $x\in Ord^{\veps}$ with
the `ordinal' coded by $x$,
and $<^{\veps}$ is denoted by $<$ when no confusion likely occurs.
Note that the ordinal $\rho_{0}$ 
is the order type of the class of ordinals in the intended model $L_{\rho_{0}}$ of $T_{1}$.
Define simultaneously 
 the classes $\calh_{\alp}(X)\subset \veps_{\rho_{0}+1}$
and the ordinals $\Psi_{\ome_{1}} (\alp)$ and $\Psi_{\rho_{0}}(\alp)$ 
for $\alp<^{\veps}\veps_{\rho_{0}+1}$ and sets $X\subset \veps_{\rho_{0}+1}$ as follows.
We see that $\calh_{\alp}(X)$ and $\Psi_{\kap} (\alp)\, (\kap\in\{\ome_{1},\rho_{0}\})$ are (first-order) definable as a fixed point in $T_{1}$.
Recall that $\mbox{Hull}(X)\subset L_{\rho_{0}}$ and
$F_{X}:\mbox{Hull}(X)\lrarw L_{\gam}$ 
for $X\subset L_{\rho_{0}}$ and a $\gam=F_{X}(\rho_{0})\leq\rho_{0}$.

\bdf\label{df:Cpsiregularsm}

$\calh_{\alp}(X)$ {\rm is defined recursively as follows.}

\benu
\item
$\{0,\ome_{1},\rho_{0}\}\cup X\subset\calh_{\alp}(X)$.

\item
 $x, y \in \calh_{\alp}(X)\Rarw x+ y,\ome^{x}\in \calh_{\alp}(X)$.

\item
$\gam\in \calh_{\alp}(X)\cap\alp
\Rarw 
\Psi_{\rho_{0}}(\gam)\in\calh_{\alp}(X)
$.

\item
$\gam\in \calh_{\alp}(X)\cap\alp
\Rarw 
x=\Psi_{\ome_{1}}(\gam)\in\calh_{\alp}(X) \spand F_{x\cup\{\ome_{1}\}}(\rho_{0})\in\calh_{\alp}(X)$.

\item

{\rm Let $A(x;y_{1},\ldots,y_{n})$ be a $\Del_{0}$-formula in the language $\{\in\}$.
For $\{\alp_{1},\ldots,\alp_{n}\}\subset\calh_{\alp}(X)$,
$\mu x.\, A(x;\alp_{1},\ldots,\alp_{n})\in\calh_{\alp}(X)$,
where $\mu x.A(x;\alp_{1},\ldots,\alp_{n})=\bet$ for 
 the least ordinal $\bet$ such that $A(\bet;\alp_{1},\ldots,\alp_{n})$
if such an ordinal exists.
Otherwise $\mu x. A(x;\alp_{1},\ldots,\alp_{n})=0$.
}

\eenu

{\rm For} $\kap\in\{\ome_{1},\rho_{0}\}$ {\rm and} $\alp<\veps_{\rho_{0}+1}$
\[
\Psi_{\kap}(\alp):=
\min\{\bet\leq\kap :  \calh_{\alp}(\bet)\cap \kap \subset \bet \}
.\]
\edf

The ordinal $\Psi_{\kap}(\alp)$ is well defined and $\Psi_{\kap}(\alp)\leq \kap$
for $\kap\in\{\ome_{1},\rho_{0}\}$.

\brem\label{rem:rec}
{\rm
In Definition \ref{df:Cpsiregularsm} let $\calh_{\alp}^{-}(X)$ denote the Skolem hull of $\{0,\ome_{1},\rho_{0}\}$
under under the functions
$+, \alp\mapsto\ome^{\alp}, \psi_{\ome_{1}}\restrict\alp, \psi_{\rho_{0}}\restrict \alp$,
with $\psi_{\kap}(\alp)=\min\{\bet\leq\kap: \calh_{\alp}^{-}(\bet)\cap\kap\subset\bet\}$ for $\kap\in\{\ome_{1},\rho_{0}\}$.
Namely $\calh_{\alp}^{-}(X)$ need not to be closed under the $\mu$-operator and the Mostowski
collapsing $F$.
Then $\psi_{\ome_{1}}(\veps_{\rho_{0}+1})$ gives the proof-theoretic ordinal of the set theory
${\sf KP}\ome+$`there exists a recursively regular ordinal', i.e.,
the theory ${\sf ID}_{2}$ of two times iterated positive inductive definitions over $\mathbb{N}$, which is
fairly weaker than the theory $T_{1}$.
}
\erem

\bprp\label{prp:definability}
Each of 
$x=\calh_{\alp}(X)$,
$y=\Psi_{\kap}(\alp)\,(\kap\in \{\ome_{1},\rho_{0}\})$
and $u=F_{z\cup\{\ome_{1}\}}(\rho_{0})=\sup\{F_{z\cup\{\ome_{1}\}}(\alp) :\alp\in\mbox{{\rm Hull}}(z\cup\{\ome_{1}\})\}$
is a $\Sig_{2}$-predicate in $\KP\ome$. 
\eprp
\bprf
Each of 
$x=\calh_{\alp}(X)$ and
$y=\Psi_{\kap}(\alp)\,(\kap\in \{\ome_{1},\rho_{0}\})$
is seen to be a $\Sig_{2}$-predicate as a fixed point so that
$\Psi_{\kap}(\alp)=x \to x<\kap$.
We see that $h=\mbox{{\rm Hull}}(z\cup\{\ome_{1}\})$ is a $Bool(\Sig_{1})$-predicate, and
the Mostowski collapsing $F^{h}$ of a set $h$ is a $\Del_{1}$-map.
\eprf

\bdf\label{df:ordterm}
{\rm Define inductively sets of ordinal terms $T(\veps_{\rho_{0}+1})$ and 
its subsets $T(\ome_{k}(\rho_{0}+1))\, (k<\ome)$
as follows. 
Each element in the sets is a term over constants $0,\ome_{1},\rho_{0}$ and function symbols
$\#,\ome, D_{0},D_{1}, F$ and $f_{A}$ for $\Del_{0}$-formulas $A$ in the language $\mathcal{L}_{2}$.
$v(a)=\alp$ designates that the value of the term $a$ is the ordinal $\alp<\veps_{\rho_{0}+1}$.
}

\benu
\item
 \benu
 \item
 {\rm
 $\{0,\ome_{1},\rho_{0}\}\subset T(\ome_{k}(\rho_{0}+1))$ for each $k<\ome$.
 }
 \item
 {\rm If $\{a_{1},\ldots,a_{n}\}\subset T(\ome_{k}(\rho_{0}+1))$ with $n>1$,
 then $a_{1}\#\cdots\#a_{n}\in T(\ome_{k}(\rho_{0}+1))$.
 }
 \item
 {\rm If $a\in T(\ome_{k}(\rho_{0}+1))$, then $\ome^{a}\in T(\ome_{k+1}(\rho_{0}+1))$.
 }
 \item
 {\rm If $a\in T(\ome_{k}(\rho_{0}+1))$, then $D_{1}(a),D_{0}(a), F(a)\in T(\ome_{k}(\rho_{0}+1))$.}
 \item
 {\rm Let $A(x;y_{1},\ldots,y_{n})\,(n\geq 0)$ be a $\Del_{0}$-formula in $\mathcal{L}_{2}$,
 and $\{a_{1},\ldots,a_{n}\}\subset T(\ome_{k}(\rho_{0}+1))$.
 Then $f_{A}(a_{1},\ldots,a_{n})\in T(\ome_{k}(\rho_{0}+1))$.
 }
 \eenu

\item
$T(\veps_{\rho_{0}+1})=\bigcup_{k<\ome}T(\ome_{k}(\rho_{0}+1))$.

\item
{\rm For $a\in T(\veps_{\rho_{0}+1})$, the \textit{value} $v(a)$ of $a$ is defined recursively as follows.
$v(0)=0$, $v(\ome_{1})=\ome_{1}$, and $v(\rho_{0})=\rho_{0}$ with ordinals in the right hand sides.
$v(a_{1}\#\cdots\#a_{n})=v(a_{1})\#\cdots\#v(a_{n})$ with the natural sum $\#$ on ordinals in the right hand side,
and $v(\ome^{a})=\ome^{v(a)}$.
$v(D_{1}(a))=\Psi_{\rho_{0}}(v(a))$, $v(D_{0}(a))=\Psi_{\ome_{1}}(v(a))$, and
$v(F(a))=F_{x\cup\{\ome_{1}\}}(\rho_{0})$ with $x=\Psi_{\ome_{1}}(v(a))$,
where $v(D_{1}(a))=y$ denotes $y<\rho_{0}\land \Psi_{\rho_{0}}(v(a),y)$ for the $\Sig_{2}$-predicate
$\Psi_{\rho_{0}}(\alp,y)\lrarw\Psi_{\rho_{0}}(\alp)=y$, and similarly for $v(D_{0}(a))=y$.
$v(f_{A}(a_{1},\ldots,a_{n}))=\mu x. A(x;\alp_{1},\ldots,\alp_{n})$ with $\alp_{i}=v(a_{i})$.
}

\eenu
\edf

Proposition \ref{prp:definability} yields the following Proposition \ref{prp:value}.
By the definition $v(D_{1}(a))=x\to x<\rho_{0}$ and
$v(D_{0}(a))=x\to x<\ome_{1}$ hold.

\bprp\label{prp:value}
For $a\in T(\veps_{\rho_{0}+1})$ and $\alp<\veps_{\rho_{0}+1}$,
$v(a)=\alp$ is a
$\Sig_{2}$-predicate in $\KP\ome$. 
\eprp

\blem\label{lem:lowerbndreg}
For {\rm each} $k<\ome$,
$T_{1}\vdash\fal a\in T(\ome_{k+1}(\rho_{0}+1))\, A(a)$ for the formula $A$ in (\ref{eq:A}).
\elem
\bprf
Let $\kap\in\{\ome_{1},\rho_{0}\}$.
By Proposition \ref{prp:definability}
 both 
$x=\calh_{\alp}(\bet)$ and
$y=\Psi_{\kap}(\alp)$ are $\Sig_{2}$-predicates.

We show that $B(\alp):\Lrarw \left(\fal a\in T(\veps_{\rho_{0}+1})(v(a)=\alp\to A(a))\right)$ is progressive.
Then $\fal\alp<\ome_{k+1}(\rho_{0}+1)\, B(\alp)$
 follow from transfinite induction up to $\ome_{k+1}(\rho_{0}+1)$, cf. Proposition \ref{prp:canonical}.\ref{prp:canonical2}.
$\fal\alp<\ome_{k+1}(\rho_{0}+1)\, B(\alp)$ yields 
$\fal a\in T(\veps_{\rho_{0}+1})\exi\alp<\ome_{k+1}(\rho_{0}+1)(v(a)=\alp)$
since $\fal a\in T(\ome_{k+1}(\rho_{0}+1))\fal\alp(v(a)=\alp\to \alp<\ome_{k+1}(\rho_{0}+1))$.
Therefore we obtain $\fal a\in T(\ome_{k+1}(\rho_{0}+1))\, A(a)$.

Assume $\fal\gam<\alp\, B(\gam)$ as our induction hypothesis.
We show $\exi x<\kap[\Psi_{\kap}(\alp)=x]$ for $\kap\in\{\ome_{1},\rho_{0}\}$.

We see from $\fal \bet<\rho_{0}\exi h[h=\mbox{{\rm Hull}}(\bet)]$, the induction hypothesis and 
$\Sig_{2}$-Collection that $\fal\bet<\rho_{0}\exi x[x=\calh_{\alp}(\bet)=\bigcup_{m}\calh_{\alp}^{m}(\bet)]$.

Define recursively ordinals $\{\bet_{m}\}_{m}$ for $\kap\in \{\ome_{1},\rho_{0}\}$ as follows.
Let $\bet_{0}=0$.
$\bet_{m+1}$ is defined to be the least ordinal $\bet_{m+1}\leq\kap$ such that
$\calh_{\alp}(\bet_{m})\cap\kap\subset\bet_{m+1}$.

We see inductively that $\bet_{m}<\kap$ using the regularity of $\kap$ and the facts that
$\fal \bet<\kap\exi x[x=\calh_{\alp}(\bet) \land card(x)<\kap]$, where
$card(x)<\kap$ designates that there exists a surjection $f:\gam\to x$ for a $\gam<\kap$ and $f\in L_{\rho_{0}}$.
Moreover $m\mapsto\bet_{m}$ is a $\Sig_{2}$-map.
Therefore $\bet=\sup_{m}\bet_{m}<\kap$ enjoys 
$\calh_{\alp}(\bet)\cap\kap\subset\bet$.
Therefore $\exi x<\kap[\Psi_{\kap}(\alp)=x]$.

Let $x=\Psi_{\ome_{1}}(\alp)<\ome_{1}$.
$F_{x\cup\{\ome_{1}\}}(\rho_{0})<\ome_{1}$ is seen from $x<cf(\ome_{1})$ and the $\Sig_{1}$-projectum
$\rho(L_{\rho_{0}})=\rho_{0}$, i.e., $\rho_{0}$ is nonprojectible, cf.\,Lemma 2.8 in \cite{liftupZF}.
\eprf
\\

In what follows let us identify the term $a\in T(\veps_{\rho_{0}+1})$ with its value $v(a)=\alp<\veps_{\rho_{0}+1}$,
and $a<b:\Lrarw v(a)<^{\veps}v(b)$.

\bdf\label{df:G}
{\rm For $a,b\in T(\veps_{\rho_{0}+1})$, a finite set $G_{a}(b)\subset T(\veps_{\rho_{0}+1})$ 
is defined recursively as follows.
\benu
\item

$G_{a}(0)=G_{a}(\ome_{1})=G_{a}(\rho_{0})=G_{a}(b)=\emptyset$ for $b<a$.

In what follows $G_{a}(b)$ is defined for $b\geq a$.
\item
$G_{a}(b_{1}\#\cdots\#b_{n})=G_{a}(f_{A}(b_{1},\ldots,b_{n}))=\bigcup\{G_{a}(b_{i}):i=1,\ldots,n\}$.
 \item
$G_{a}(\ome^{b})=G_{a}(b)$.
 \item
 $G_{a}(D_{i}(b))=G_{a}(F(b))=\{b\}\cup G_{a}(b)$.
\eenu
}
\edf

\bprp\label{prp:G}
For $a,b\in T(\veps_{\rho_{0}+1})$,
$G_{a}(b)<c\Rarw b\in\calh_{c}(a)$.
\eprp

\section{Finite proof figures}\label{sect:finiteproof}

In this section  an extension $T_{c}(\ome_{1})$ of the theory $T^{ord}(\ome_{1})$ with individual constants and function constants
is formulated in one-sided sequent calculus, and
an ordinal assignment to sequents occurring in proofs are defined in subsection \ref{subsec:ordinalassignment}.

In this section \ref{sect:finiteproof} and the next section \ref{sect:consisprf}
we work in the theory 
$\KP\ome+(V=L)+\left(\fal a\in T(\veps_{\rho_{0}+1})\, A(a)\right)$.

The language $\calL_{c}$ of $T_{c}(\ome_{1})$ is obtained from $\mathcal{L}_{2}$
by adding names (individual constants) $c_{t}$
of each term $a\in T(\veps_{\rho_{0}+1})$ with $v(a)<\rho_{0}$.
The constant $c_{a}$ is identified with the term $a\in T(\veps_{\rho_{0}+1})$.
Formulas are assumed to be in negation normal form.

\bdf
\benu
\item
{\rm A \textit{literal} is one of atomic formulas $s<t, R^{\cala}(s,t)$, $P(t_{0},t_{1})$,
$P_{\rho_{0}}(t)$ 
or their negations.}

\item
{\rm The \textit{truth}  of closed literals is defined as follows.}
$s< t$ {\rm is} \textit{true} {\rm if} $v(s)< v(t)$. 
$R^{\cala}(s,t)$ {\rm is} \textit{true} {\rm if} $R^{\cala}(v(s),v(t))$ {\rm holds.}
$P(t_{0},t_{1})$ {\rm is} \textit{true} {\rm if} $v(t_{0})=x=\Psi_{\ome_{1}}(\bet)$ {\rm and}
$v(t_{1})=F_{x\cup\{\ome_{1}\}}(\rho_{0})$ {\rm for some} $\bet$.
$P_{\rho_{0}}(t)$ {\rm is} \textit{true} {\rm if} $v(t)=\Psi_{\rho_{0}}(\bet)$ {\rm for some $\bet$.}
{\rm A closed literal $\lnot L$ is \textit{true} if $L$ is not true.}

\item
{\rm An} $E$\textit{-formula} {\rm is either a literal or a formula of one of the shapes $A_{0}\lor A_{1}, \exi x\, A(x)$.}
\eenu
\edf
By \textit{$\Del_{0}$-formula} we mean a bounded formula in the language $\mathcal{L}_{c}$ in which predicates $P,P_{\rho_{0}}$
do not occur.
The \textit{truth} of $\Del_{0}$-sentences is defined from one of literals.
A \textit{$\Sig_{1}$-formula} or a \textit{$\Pi_{1}$-formula} is defined similarly.
These formulas are obtained from formulas in $\mathcal{L}_{1}\cup\{\ome_{1}\}$ by substituting $\mathcal{L}_{c}$-terms for variables.

\textit{Proof figures} are constructed from the following axioms and inference rules in $T_{c}(\ome_{1})$.
Relations between occurrences $A,B$ of formulas in a proof such as `$A$ is a \textit{descendant} of $B$'
or equivalently `$B$ is an \textit{ancestor} of $A$',
and `an occurrence of inference rule is \textit{implicit} or \textit{explicit}'
 are defined as in \cite{ptMahlo}.
\\

\noindent
[{\bf Axioms}]
\[
\infer[(ax)]{\Gam,A}{}
\]
where $A$ is either a true closed literal or a true closed $\Del_{0}$-formula or
a $\Del_{0}$-axiom whose universal closure is an axiom for the constants
$0,<,+,\cdot,\lam x.\ome^{x}$.

\[
\infer[(taut)]{\Gam,\lnot A, A}{}
\msfiv\mbox{   for literals and $\Del_{0}$-formulas $A$.}
\]
This means $ \dg(A)=1$ in Definition \ref{df:depthfml} below.
If there occurs no fee variable in an axiom $(ax), (taut)$,
then it contains either a true literal or a true $\Del_{0}$-sentence.
\[
\infer[(P\exi)]{\Gam, (s\not<\ome_{1},) \exi x, y<\ome_{1}[s<x\land P(x,y)]}{}
\]
Cf.\,(\ref{eq:Z2ord}). When $s<\ome_{1}$ is a true literal, $s\not<\ome_{1}$ may be absent.

Cf.\,(\ref{eq:Z5ord}).
\[
\infer[(P_{\rho_{0}}\exi)]{\Gam,  \exi x[s<x\land P_{\rho_{0}}(x)]}{}
\]
\\

\noindent
[{\bf Inference rules}]
In each case the main (principal) formula is assumed to be in the lower sequent $\Gam$.
Namely
$(A_{0}\lor A_{1})\in\Gam$ in $(\lor)$, $(A_{0}\land A_{1})\in\Gam$ in $(\land)$,
$(\exi x\, A(x))\in\Gam$ in $(\exi)$, $(\exi x<t\, A(x))\in\Gam$ in $(b\exi)$, $(\fal x\, A(x))\in\Gam$ in $(\fal)$, $(\fal x<t\,A(x))\in\Gam$ in $(b\fal)$.

The variable $x$ in $(\fal),(b\fal)$ is an eigenvariable.

\[
\infer[(\lor)]{\Gam}{\Gam,A_{i}}
\:
\infer[(\land)]{\Gam}{\Gam, A_{0} & \Gam,A_{1}}
\]
\[
\infer[(\exi)]{\Gam}{\Gam, A(s)}
\:
\infer[(b\exi)]{(s\not<t,)\Gam}
{
\Gam,A(s)
}
\:
\infer[(\fal)]{\Gam}{\Gam, A(x)}
\:
\infer[(b\fal)]{\Gam}{\Gam,x\not<t, A(x)}
\]
where in $(b\exi)$, the formula $s\not <t$ may be absent when $s<t$ is a closed true literal.
\[
\infer[(ind)]{(s\not< t,)\Gam}{\Gam,\lnot\fal x< y A(x), A(y) & \Gam,\lnot A(s)}
\]
where 
$s\not< t$ may be absent in the lower sequent when
$s< t$ is a true closed formula.
The formula $A(x)$ is the \textit{induction formula}, and
the term $t$ is the \textit{induction term} of the $(ind)$.
The variable $y$ is the \textit{eigenvariable} of the rule $(ind)$.

\[
\infer[(cut)]{\Gam,\Lam}{\Gam,\lnot A & A,\Lam}
\]
$A$ is an $E$-formula called the \textit{cut formula} of the $(cut)$.
\[
\infer[(Rfl)]{\Gam}
{
\Gam, \fal x< t\, A(x)
&
t\not< y,\exi x< t\, \lnot A^{(y)}(x),\Gam
}
\]
$t$ is a term, $y$ is an eigenvariable,  and 
$A(x)\equiv(\exi z\exi w[P_{\rho_{0}}(z)\land B(x)])\, (B\in\Del_{0})$, 
$A^{(y)}(x):\equiv(\exi z<y\exi w<y[P_{\rho_{0}}(z)\land B(x)])$,
cf.\,(\ref{eq:Z6ord}).

\[
\infer[(P\Sig_{1})]{\Gam, (\lnot P(t_{0},t_{1}), s\not<t_{0},)  \vphi^{t_{1}}[t_{0},s]}{\Gam,\vphi[\ome_{1},s]}
\]
$\vphi$ is an arbitrary $\Sig_{1}$-formula in the language $\mathcal{L}_{1}$ with predicates $R^{\cala}$,
cf.\,(\ref{eq:Z1ord}).
When $P(t_{0},t_{1})$ or $s<t_{0}$ is a true literal, their negations $\lnot P(t_{0},t_{1}), s\not<t_{0}$ may be absent.

\[
\infer[(P_{\rho_{0}}\Sig_{1})]{\Gam,(\lnot P_{\rho_{0}}(t), s\not<t,) \vphi^{t}[s]}{\Gam,\vphi[s]}
\]
$\vphi$ is an arbitrary $\Sig_{1}$-formula in the language $\mathcal{L}_{1}$, cf.\,(\ref{eq:Z4ord}).
When $P_{\rho_{0}}(t)$ or $s<t$ is a true literal, their negations $\lnot P_{\rho_{0}}(t), s\not<t$ may be absent.

\[
\infer[(h)]{\Gam,\Del}{\Gam}
\]

Let $c\oplus \alp:=c\#\alp$.
$\oplus$ is used as a punctuation mark.

\[
\infer[(D_{1})_{\alp}]{\Lam,\Gam^{(\alp)}}{\Lam, \Gam}
\]
{\rm where} 
$\alp=D_{1}(c_{1}\oplus\alp_{1})$ for some $c_{1}\oplus\alp_{1}$, and $\Gam^{(\alp)}=\{C^{(\alp)}:C\in\Gam\}$.
Each formula in $\Gam$ is one of the closed formulas 
$\fal x<t\, A(x)$, $A(s_{0})$, and $\exi w[P_{\rho_{0}}(s_{1})\land B(s_{0},s_{1},w)]$,
where $B$ is a $\Del_{0}$-formula, 
$A(x)\equiv(\exi z\exi w[P_{\rho_{0}}(z)\land B(x,z,w)])$
$t,s_{0},s_{1}$ are closed terms.
Each \textit{implicit} formula in $\Lam$ is a bounded sentence.
Note that there occurs no unbounded universal quantifier in any implicit formula in
$\Lam\cup\Gam$.

\[
\infer[(D_{0})_{\alp}]{\Lam}{\Lam}
\]
where each formula in $\Lam$ is either a false closed $\Del_{0}$-formula or
a closed subformula of a $\Sig_{2}$-sentence $\exi x \fal y \, B(x,y)$.
$\alp=D_{0}(c_{0}\oplus \alp_{0})$ for some $c_{0}\oplus \alp_{0}$.

\subsection{Ordinal assignment}\label{subsec:ordinalassignment}
In this subsection let us define ordinal assignments.

\bdf\label{df:height}
{\rm The \textit{height} $h(\Gam)=h(\Gam;\calP)<\ome\cdot 2$ of sequents $\Gam$ in a proof figure $\calP$.}
\benu
\item
$h(\Gam)=0$ {\rm if $\Gam$ is the end-sequent of $\calP$.}

\item
$h(\Gam)=\ome\cdot i$ {\rm if $\Gam$ is the upper sequent of a $(D_{i})$ with $i=0,1$.}

\item
$h(\Gam)=h(\Del)+1$ {\rm if $\Gam$ is the upper sequent of an $(h)$ with its lower sequent $\Del$.}

\item
$h(\Gam)=h(\Del)$ {\rm if $\Gam$ is an upper sequent of a rule other than $(h)$ and $(D_{i})$ with its lower sequent $\Del$.}
\eenu
{\rm Let
$h_{0}(\Gam)=h(\Gam)$ if $h(\Gam)<\ome$.
$h_{0}(\Gam)=h(\Gam)-\ome$ if $h(\Gam)\geq \ome$.}
\edf

\bdf\label{df:depthfml}{\rm The \textit{degree} $ \dg(A)<\ome$ of formulas $A$.}
\benu
\item
$ \dg(A)=1$ {\rm if $A$ is either a literal or a $\Del_{0}$-formula.}

{\rm In what follows} $A$ {\rm is neither a literal nor a $\Del_{0}$-formula.}

\item
$ \dg(A)= \dg(A_{0})+ \dg(A_{1})+2$ {\rm if $A\equiv(A_{0}\lor A_{1}), (A_{0}\land A_{1})$.}

\item
$ \dg(A)= \dg(B)+2$ {\rm if $A\equiv(\exi x\, B(x)),(\fal x\, B(x))$.}

\item
$ \dg(A)= \dg(B)+2$ {\rm if $A\equiv(\exi x< t\, B(x)),(\fal x< t\, B(x))$.}
\eenu
\edf

\bdf
{\rm A proof figure is said to be \textit{height regulated} if it enjoys the following conditions:}
\bdes

\item[(h1)] 
{\rm There occurs no free variable in any sequent $\Gam$ if $h(\Gam)<\ome$.}

\item[(h2)] 

{\rm Let $\Gam,\exi x[s<x\land P_{\rho_{0}}(x)]$ be an axiom $(P_{\rho_{0}}\exi)$ in $\mathcal{P}$,
and $J$ be a $(cut)$ whose cut formula is a descendant $C\equiv(\exi x[s<x\land P_{\rho_{0}}(x)])$
of $C$ in the axiom.
Then $h(\Del)\geq\ome$ for the upper sequent $\Del$ of the $(cut)\, J$.}

\item[(h3)] 
{\rm For any $(cut)$ in $\calP$,
$ \dg(C)\leq h_{0}(\Gam,\Del)$ for its cut formula $C$ and the lower sequent $\Gam,\Del$.}

\item[(h4)]
{\rm For any $(ind)$ in $\calP$
\[
\infer[(ind)]{(s\not< t,)\Gam}
{
\Gam,\lnot\fal x< y A(x), A(y) 
& 
\Gam,\lnot A(s)
}
\]
$\ome+ \dg(\fal x< s A(x))\leq h(s\not< t, \Gam)$ holds, and
there are no nested $(ind)$ rules, i.e., there occurs no $(ind)$ above the rule $(ind)$.}

\item[(h5)]
{\rm There exists a rule $(D_{1})$ below a $(Rfl)$. 
Let $J$ be the lowest such rule $(D_{1})$ with the lower sequent $\Del$.
Then $h(\Del)\geq \dg(\exi x< t \lnot A^{(y)}(x))$.}
\[
\infer[(D_{1})\, J]{\Del}
{
 \infer*{\cdots}
   {
    \infer[(Rfl)]{\Gam}
     {
     \Gam, \fal x< t\, A(x)
     &
    t\not< y, \exi x< t \lnot A^{(y)}(x),\Gam
     }
 }
}
\]

\item[(h6)]
{\rm If a rule $(D_{1})\, J_{0}$ is above another $(D_{1})\, J_{1}$, then the only rules between $J_{0}$ and $J_{1}$ are
$(D_{1})$'s.}

\item[(h7)]
{\rm $\calP$ ends with an inference rule $(D_{0})$.}

\edes
\edf
Let $\calP$ be a height regulated proof with a rule $(D_{1})$.
By {\bf (h6)}, rules $(D_{1})$ occur consecutively.
\[
\infer[(D_{1})]{\Gam_{0}}
{
\infer{\Gam_{1}}
 {
  \infer*{}
 {
 \infer[(D_{1})]{\Gam_{n-1}}{\Gam_{n}}
 }
}
}
\]
with $h(\Gam_{0})<\ome$.

Let us assign an ordinal term $c_{1}\in T(\veps_{\rho_{0}+1})$ to each lowest rule $(D_{1})$ occurring in $\calP$.
$c_{1}$ is the \textit{stock} of each rule $(D_{1})$ in the consecutive series.
Also a term $c_{0}$ is assigned to the last rule $(D_{0})$, the \textit{stock} of the $(D_{0})$, cf.\,{\bf (h7)}.
Such an assignment ${\sf c}$ is said to be a \textit{stock assignment} for $\calP$.


\bdf\label{df:ordinalassignment}
{\rm Given a stock assignment ${\sf c}$, we assign an ordinal term $o(\Gam)=o(\Gam;\calP,{\sf c})\in T(\veps_{\rho_{0}+1})$
 to each occurrence of a sequent $\Gam$ in a proof figure $\calP$.
Let us write 
\[
\infer{\Gam; b}{\cdots\Gam_{i}; a_{i}\cdots}
\]
when the lower sequent $\Gam$ receives an ordinal term $b$, i.e., $o(\Gam)=b$, and
$o(\Gam_{i})=a_{i}$ for upper sequents $\Gam_{i}$.}
\\

\noindent
{\bf Axioms} 
{\rm 
If $\Gam$ is one of axioms $(ax), (taut),(P\exi), (P_{\rho_{0}}\exi)$, then $o(\Gam)=1=\ome^{0}$.
}
\\

\noindent
{\bf Rules}.
{\rm Let $\Gam$ be the lower sequent of a rule $J$ with its upper sequents $\Gam_{i}$:}
\[
\infer[J]{\Gam}{\cdots\Gam_{i}\cdots}
\]

\benu

\item
{\rm $J$ is one of the rules $(P\Sig_{1})$ or $(P_{\rho_{0}}\Sig_{1})$:}
 $o(\Gam)=o(\Gam_{0})$.
 
\item\label{df:ordinalassignmentbfal}
{\rm $J$ is one of rules $(\lor), (b\exi), (\exi), (b\fal), (\fal)$:
$o(\Gam)=o(\Gam_{0})+1$.}

\item
{\rm $J$ is one of rules $(\land), (cut),(Rfl)$:
$o(\Gam)=o(\Gam_{0})\# o(\Gam_{1})$.}


\item
{\rm $J$ is an $(h)$:
 $o(\Gam)=\ome^{o(\Gam_{0})}$. 
It is convenient for us to write
$D_{2}(0\oplus\alp):=\ome^{\alp}$ and $(D_{2}):=(h)$.
}


\item
{\rm $J$ is an $(ind)$:
\[
\infer[(ind)]{(s\not< t,) \Gam;b}{\Gam,\lnot\fal x< y A(x), A(y);a_{0} & \Gam,\lnot A(s);a_{1}}
\]
Let $mj(t)=\rho_{0}$ if $t$ is not closed.
Otherwise $mj(t)=t$.
Then
$b= (a_{0}+a_{1}+2)\times mj(t)$ for  the natural product $\times$, cf.\,{\bf (p1)} below.}

\item
{\rm $J$ is a rule $(D_{1})$:}
\[
o(\Gam)=
\left\{
\begin{array}{ll}
D_{1}(c_{1}\oplus \ome^{o(\Gam_{0})}) & \mbox{{\rm if }} h(\Gam)<\ome
\\
o(\Gam_{0}) & \mbox{{\rm if }}  h(\Gam)=\ome
\end{array}
\right.
\]
{\rm where $c_{1}$ is the stock ${\sf c}(J)$ of the rule $J$ under the stock assignment.}

\item
{\rm $J$ is a rule $(D_{0})$:}
\[
o(\Gam)=D_{0}(c_{0}\oplus o(\Gam_{0})) 
\]
{\rm where $c_{0}$ is the stock ${\sf c}(J)$ of the rule $J$ under the stock assignment.}

\eenu
{\rm Finally let
$o(\calP)=o(\Gam_{end})$ for the end-sequent $\Gam_{end}$ of $\calP$.}
\edf


\blem\label{lem:depth}{\rm (Tautology lemma)}\\
For any formula $A(x)$, any $\Gam$ and any term $t$,
there exists a proof $\calP$ of $\Gam,\lnot A(t),A(t)$ such that 
$o(\Gam,\lnot A(t),A(t);\calP,{\sf c})= \dg(A(x))$ for any stock assignment ${\sf c}$.
\elem
\bprf 
We see the assertion by induction on $ \dg(A)$.
\eprf

\bprp\label{prp:oa}
Let $\Gam$ be a sequent in a proof $\calP$ with a stock assignment ${\sf c}$, and 
$b\in G_{a}(o(\Gam;\calP,{\sf c}))$.
\benu
\item\label{prp:oa.1}
If $h(\Gam)\geq\ome$, then there exists a closed induction term $t$ occurring above $\Gam$
such that $b\in G_{a}(t)$.
\item\label{prp:oa.2}
Let $h(\Gam)<\ome$ for the sequent $\Gam$ other than the end-sequent.
Then either there exists a closed induction term $t$ occurring above $\Gam$
such that $b\in G_{a}(t)$, or there exists a lowest rule $(D_{1})\, J$ such that 
$b=c\oplus \ome^{o(\Del)}$ for $c={\sf c}(J)$ and the upper sequent $\Del$ of $J$.
\eenu
\eprp

\bdf\label{df:proofoa}
{\rm A proof figure $\calP$ together with a stock assignment ${\sf c}$
is a \textit{proof with stock} 
if the following conditions are met.}
\bdes
\item[(p0)] 
{\rm $\calP$ is height regulated.}

\item[(p1)]
{\rm For any $(ind)$ occurring in $\calP$}
\[
\infer[(ind)]{(s\not< t,) \Gam;a}{\Gam,\lnot\fal x< y A(x), A(y);a_{0} & \Gam,\lnot A(s);a_{1}}
\]
$ \dg(A(y))=a_{1}$ {\rm and} $a_{0}<\ome$.

\item[(p2)]
{\rm Let $J$ be a rule $(D_{i})_{\alp}$ with an ordinal $\alp=D_{i}(\alp_{0})$
and its stock $c={\sf c}(J)$ occurring in $\calP$.
Then $\fal d\in T(\veps_{\rho_{0}+1})[ G_{D_{i}(c\oplus d)}(c)<c]$,
$\alp_{0}\geq c\oplus \ome^{o(\Gam)}$ and $\alp\geq D_{i}(c\oplus \ome^{o(\Gam)})$ with
the upper sequent $\Gam$ of $J$.}
 \bdes
  
 \item[(p2.1)]
{\rm Let $t$ be a closed term occurring above the rule $(D_{i})\, J$.
Then $\fal d\in T(\veps_{\rho_{0}+1})[ G_{D_{i}(c\oplus d)}(t)<c]$,
where by a closed term occurring in a proof we mean to include a closed subterm in a term occurring in the proof.}

\item[(p2.2)]
{\rm
Let $i=0$ and $J_{1}$ be a rule $(D_{1})_{\bet}$ occurring above the rule $J$ 
with an ordinal $\bet=D_{1}(\bet_{0})$ and its stock $c_{1}$.
Then $\bet_{0}<c$ and
$\fal d\in T(\veps_{\rho_{0}+1})[ G_{D_{0}(c\oplus d)}(c_{1})<c]$ for the stock $c_{1}={\sf c}(J_{1})$ of the rule
$(D_{1})\, J_{1}$.
}



 \edes

\edes
\edf

\brem
{\rm
The condition $\fal d\in T(\veps_{\rho_{0}+1})[ G_{D_{i}(c\oplus d)}(t)<c]$ in {\bf (p2.1)}
yields $t\in\bigcap\{H_{c\oplus d}(D_{i}(c\oplus d)): d\in T(\veps_{\rho_{0}+1})\}$ by Proposition \ref{prp:G}.
Also the condition is equivalent to $G_{D_{i}(c)}(t)<c$ under $(\Pi_{1}\mbox{{\rm -Collection}})$
since $(\Pi_{1}\mbox{{\rm -Collection}})$ yields $\calh_{c}(D_{i}(c\oplus d))\subset\calh_{c\oplus d}(D_{i}(c\oplus d))$,
and
$D_{i}(c)\leq D_{i}(c\oplus d)$.
}
\erem

\blem\label{lem:delta0elim}{\rm (False literal elimination)}\\
Let $A$ be a false closed literal, and
$\calP$ a proof of $\Gam,A$.
Then there exists a proof $\calP'$ of $\Gam$ such that
$o(\Gam;\calP',{\sf c})=o(\Gam,A;\calP,{\sf c})$ for any stock assignment ${\sf c}$.
\elem
\bprf
Eliminate the ancestors $A$ of $A$ to get a proof $\calP'$ of $\Gam$.
Consider a $(P\Sig_{1})$.
\[
\infer[(P\Sig_{1})]{\Gam, (\lnot P(t_{0},t_{1}), s\not<t_{0},)  \vphi^{t_{1}}[t_{0},s];a}{\Gam,\vphi[\ome_{1},s];a}
\]
If one of literals $\lnot P(t_{0},t_{1}), s\not<t_{0}$ is a false ancestor of $A$, then eliminate it from the lower sequent.
The case $(P_{\rho_{0}}\Sig_{1})$ is similar.
\eprf
\\

The following Lemma \ref{lem:reflect} yields Theorem \ref{th:main}.

\blem\label{lem:reflect}
Let $\calP$ be a proof with stock.
Then $\bigvee\Gam$ is true for the end-sequent $\Gam$ of $\calP$.
\elem

Lemma \ref{lem:reflect} is shown by induction on ordinals $o(\calP)<\ome_{1}$ 
in section \ref{sect:consisprf}.

\subsection{Initial proofs}\label{subsec:embedding}

\blem\label{lem:embed}
Suppose that $T^{ord}(\ome_{1})$ proves a $\Sig_{2}$-sentence $\exi x\fal y\, C_{0}(x,y)$ 
with a $\Del_{0}$-formula $C$ in the language $\mathcal{L}_{1}$.
Then there exists a proof $\calP_{0}$ of the $\Sig_{2}$-sentence
$\exi x\fal y\, C_{0}(x,y)$ with a stock assignment ${\sf c}_{0}$
such that $(\calP_{0},{\sf c}_{0})$ is a proof with stock.
\elem

Suppose that $T^{ord}(\ome_{1})$ proves a $\Sig_{2}$-sentence $\exi x\fal y\, C_{0}(x,y)$. 
We show that there exists a proof $\calP_{0}$ of $\exi x\fal y\, C_{0}(x,y)$ and a stock assignment ${\sf c}_{0}$
such that $(\calP_{0},{\sf c}_{0})$ is a proof with stock.

Let $\mathcal{Q}_{0}$ be a proof figure of $\exi x\fal y\, C_{0}(x,y)$ from 
axioms (\ref{eq:Z1ord}), (\ref{eq:Z2ord}), (\ref{eq:Z4ord}), (\ref{eq:Z5ord}) and (\ref{eq:Z6ord}),
and axioms in $T_{2}$ other than $\Pi_{2}$-Reflection.


Each leaf in $\mathcal{Q}_{0}$ is either a logical one $(taut)$ or 
one of axioms (\ref{eq:Z1ord}), (\ref{eq:Z2ord}), (\ref{eq:Z4ord}), (\ref{eq:Z5ord}) and (\ref{eq:Z6ord}),
and axioms in $T_{2}$ other than $\Pi_{2}$-Reflection.
Inference rules in $\mathcal{Q}_{0}$ are logical ones, $(\lor),(\land),(\exi),(\fal)$ 
and $(cut)$.
Let us depict pieces of proofs of each leaf in $\mathcal{Q}_{0}$ except $(taut)$'s.

\benu
\item
When $\fal\vec{x}\, A$ is the universal closure of an axiom  in $T_{2}$
except Foundation and $\Pi_{2}$-Reflection schema, replace the leaf $\Gam,\fal\vec{x}\, A$ by 
\[
 \infer[(\fal)]{\Gam,\fal\vec{x}\, A}
 {
  \infer[(ax)]{\Gam,A; 1}{}
 }
\]

\item
Leaves for axioms (\ref{eq:Z1ord}), (\ref{eq:Z2ord}), (\ref{eq:Z4ord}) and (\ref{eq:Z5ord}) are derived from inference rules 
$(P\Sig_{1})$, $(P\exi)$, $(P_{\rho_{0}}\Sig_{1})$ and $(P_{\rho_{0}}\exi)$, resp.
\[
\infer[(\lor),(\fal)]{\fal x,y,a(\lnot P(x,y)\lor a\not<x\lor\lnot\vphi[\ome_{1},a]\lor\vphi^{y}[x,a]); 12}
{
 \infer[(P\Sig_{1})]{\lnot P(x,y), a\not<x, \lnot\vphi[\ome_{1},a], \vphi^{y}[x,a]; 3 }
 {
  \infer*{\lnot P(x,y), a\not<x, \lnot\vphi[\ome_{1},a], \vphi[\ome_{1},a];3 }{}
 }
}
\]
with 6 times $(\lor)$, 3 times $(\fal)$, and $\dg(\vphi)=3$.
\[
\infer[(b\fal)]{\fal a<\ome_{1}(\exi x, y<\ome_{1}[a<x\land P(x,y)]); 2}
{
 \infer[(P\exi)]{a\not<\ome_{1},\exi x,y<\ome_{1}[a<x\land P(x,y)]; 1}{}
}
\]
where the formula $\fal a<\ome_{1}(\exi x, y<\ome_{1}[a<x\land P(x,y)])$ is not a $\Del_{0}$-formula.

\[
\infer[(\lor),(\fal)]{\Gam,\fal x,y(\lnot P_{\rho_{0}}(x) \lor y\not<x\lor  \lnot\vphi[y]\lor \vphi^{x}[y]);11}
{
 \infer[(P_{\rho_{0}}\Sig_{1})]{\Gam,\lnot P_{\rho_{0}}(x), y\not<x, \lnot\vphi[y], \vphi^{x}[y]; 3}
 {
 \infer*{\lnot\vphi[y], \vphi[y]; 3}{}
 }
}
\]
with 6 times $(\lor)$, 2 times $(\fal)$, and $\dg(\vphi)=3$.

\[
\infer[ (\fal)]{\Gam,\fal y(\exi x[y<x \land P_{\rho_{0}}(x)]); 2}
{
 \infer[(P_{\rho_{0}}\exi)]{\Gam, \exi x[y<x\land P_{\rho_{0}}(x)]; 1}{}
}
\]
\item
Leaves for transfinite induction schema are replaced by

{\scriptsize
\[
\hspace{-10mm}
\infer[(\fal),(\lor)]{\Gam,\fal y(\fal x< y\, A(x)\to A(y))\to\fal y\, A(y); d_{1}\times\rho_{0}+d_{0}+3}
{
\infer[(cut)]{\Gam,\lnot Prg, A(y); d_{1}\times\rho_{0}+d_{0}+1}
{
 \infer[(b\fal)]{\Gam,\lnot Prg,\fal x< y\, A(y)}
  {
  \infer[(ind)]{x\not< y, \Gam,\Del; d_{1}\times\rho_{0}}
   {
    \infer[(\land),(\exi)]{\Gam,\lnot Prg,\lnot\fal x< y A(x), A(y); d_{0}}
    {
    \infer*{\Gam,\fal x< y\, A(x),\lnot\fal x< y A(x);d}{}
    &
    \infer*{\Gam,\lnot A(y),A(y) ; d^{\prime}}{}
    }
    & 
    \hspace{-2mm}
   \infer*{\Gam,\Del,A(x),\lnot A(x);d^{\prime}}{}
    }
   }
  &
  \hspace{-25mm}
  \infer*{\Gam,\lnot Prg,\lnot\fal x< y\, A(x),A(y); d_{0}}{}
 }
}
\]
}
where $\Del=\{\lnot Prg, A(x)\}$ with $Prg\equiv(\fal y(\fal x< y\, A(x)\to A(y)))$ and 
$d= \dg(\fal x< y\, A(x)), d^{\prime}=\dg(A(x))=\max\{d-1,1\}$, $d_{0}=d+d^{\prime}+1$, and
$d_{1}=d_{0}+d^{\prime}+2$.
Also $\rho_{0}=mj(y)$.

Observe that this piece enjoys the condition {\bf (p1)}.

\item
Leaves for (\ref{eq:Z6ord}) are replaced by

\[
\infer[(\lor),(\fal)]{\Gam,\fal z[\fal x< z\, A(x) \to \exi y\fal x< z\,\lnot A^{(y)}(x)]; 24}
{
 \infer[(Rfl)]{\Gam,\Del; 21}
 {
  \infer*{\Gam, \Del,\fal x< z\, A(x);10}{}
 &
  \infer[(\exi)]{z\not< y, \exi x< z\, \lnot A^{(y)}(x),\Gam,\Del; 11}
  {
   \infer*{\exi x< z\, \lnot A^{(y)}(x),\fal x< z\, A^{(y)}(x),\Gam; 10}{}
   }
 }
}
\]
for $\Del=\{\lnot \fal x< z\, A(x),\exi y\fal x< z\, A^{(y)}(x)\}$ and 
$10=\dg(\fal x< z\, A(x))=\dg(\fal x< z\, A^{(y)}(x))$.

\eenu

Finally consider a $(cut)$:
\[
\infer[(cut)]{\Gam,\Del}
{
\Gam,\lnot A
&
A,\Del
}
\]
Replace it by
\[
\infer[(cut)]{\Gam,\Del;a_{0}\# a_{1}}
{
\Gam,\lnot A;a_{0}
&
A,\Del;a_{1}
}
\]
Let $\mathcal{Q}_{1}$ be the proof obtained from $\mathcal{Q}_{0}$
as described above with an ordinal $b$ constructed from $1, n\times\rho_{0}$ and $\#$.
Note that there occurs no inference rules $(D_{i})$ for $i=0,1$ in the constructed $\mathcal{Q}_{1}$,
and $G_{a}(b)=G_{a}(t)=\emptyset$ 
for any $a$ and any closed term $t$ occurring in $\mathcal{Q}_{1}$, cf.\,{\bf (p2.1)}.

Let $k\geq 10$
be a positive integer such that $k\geq \dg(C)$ for any cut formula $C$ occurring in $\mathcal{Q}_{1}$,
$k\geq\dg(\fal x<y\, A(x))$ for any induction formula $A(x)$ occurring in $\mathcal{Q}_{1}$.
Then add $k$-times $(h)$'s to get a proof $\mathcal{Q}_{2}$:

\[
\mathcal{Q}_{2}=
\left.
\begin{array}{c} 
\infer[(h)]{\exi x\fal y\, C_{0}(x,y);b_{1}}
{
  \infer*[\mathcal{Q}_{1}]{\exi x\fal y\, C_{0}(x,y); b}{}
 }
\end{array}
\right.
\]
where $b_{1}=\ome_{k}(b)$ with the number $k$ of $(h)$'s.
The conditions {\bf (h2)}, {\bf (h3)} and {\bf (h4)} are fulfilled with the proof $\mathcal{Q}_{2}$.

Finally let
\[
\mathcal{P}_{0}=
\left.
\begin{array}{c} 
 \infer[(D_{0})_{\alp_{0}}]{\exi x\fal y\, C_{0}(x,y);\alp_{0}}
 {
 \infer[(h)]{\exi x\fal y\, C_{0}(x,y);b_{0}}
 {
  \infer[(D_{1})_{\alp_{1}}]{\exi x\fal y\, C_{0}(x,y);\alp_{1}}
  {
    \infer*[\mathcal{Q}_{2}]{\exi x\fal y\, C_{0}(x,y);b_{1}}{}
    }
  }
 }
\end{array}
\right.
\]
where $\alp_{1}=D_{1}(0\oplus \ome^{b_{1}})$ with the empty stock $0$,
and another $k$-times $(h)$'s are attached below the $(D_{1})_{\alp_{1}}$.
The conditions {\bf (h5)} and {\bf (h6)} are fulfilled with the introduced rule $(D_{1})_{\alp_{1}}$.
For {\bf (h5)} note that $k\geq 10=\dg(\fal x<z\, A^{(y)}(x))$ for the formula 
$A^{(y)}(x)\equiv(\exi z<y[P_{\rho_{0}}(z)\land \exi w<y\, B(x)])\, (B\in\Del_{0})$
 in the inference rule $(Rfl)$.
$b_{0}=\ome_{k}(\alp_{1})$
and $\alp_{0}=D_{0}(c_{0}\oplus b_{0})$ with $c_{0}=\ome^{b_{1}}+1$.
Then $G_{a}(c_{0})=\emptyset$ for any $a$, and the condition {\bf (p2)} is enjoyed for $\mathcal{P}_{0}$.

$\calP_{0}$ with the stocks $0,c_{0}$ is a proof with stock defined in Definition \ref{df:proofoa}.
Since there occurs no constant other than $0,\ome_{1}$ in $\calP_{0}$, the condition {\bf (p2.1)} holds vacuously.
This shows Lemma \ref{lem:embed}.

\section{Reductions on finite proof figures}\label{sect:consisprf}

In what follows let $(\calP,{\sf c})$ be a proof with a stock assignment ${\sf c}$.
Let $\Gam_{end}$ be the end-sequent of $\calP$.
Assuming that $\bigvee\Gam_{end}$ is false,
we rewrite $(\calP,{\sf c})$ to another proof $(\calP',{\sf c}')$ with stock
so that $o(\calP')<o(\calP)$ and $\bigvee\Gam_{end}'$ is false for the end-sequent $\Gam_{end}'$ of $\calP'$.
This proves Lemma \ref{lem:reflect}.

In each case below the new stock assignment ${\sf c}'$ for the new proof $\calP'$ 
is defined obviously from the ${\sf c}$ otherwise stated.

\bdf\label{df:mbranchtop}
{\rm The \textit{main branch} of a proof figure $\calP$ is a series $\{\Gam_{i}\}_{i\leq m}$ of occurrences of sequents in $\calP$ such that:}
\benu
\item
{\rm $\Gam_{0}$ is the end-sequent of $\calP$.}
\item
{\rm For each $i<m$, $\Gam_{i+1}$ is the \textit{rightmost} upper sequent of a rule $J_{i}$ with its lower sequent $\Gam_{i}$, and
$J_{i}$ is one of the rules $(cut), (h)$,
and $(P\Sig_{1}), (P_{\rho_{0}}\Sig_{1}), (D_{i})\, (i=0,1)$.}

\item
{\rm $\Gam_{m}$ is either an axiom or the lower sequent of one of rules 
$(\lor),(\land),(\exi), (\fal)$,
$(b\exi),(b\fal),(ind),(Rfl)$.}

\eenu
{\rm $\Gam_{m}$ is said to be the \textit{top} of the main branch of $\calP$.}
\edf

Let $\Phi$ denote the top of the main branch of the proof $\calP$ with stock assignment ${\sf c}$. 
Observe that we can assume that $\Phi$ contains no free variable.

\subsection{top=axiom}\label{subsec:topaxiom}
In this subsection we consider the cases when the top $\Phi$ is an axiom.
\\

\noindent
{\bf Case 1}.
The top $\Phi=A,\Del_{0}$ is either an $(ax)$ or a $(taut)$.
Then $\Phi$ contains a true $\Del_{0}$-formula $A$ or a true literal $A=(\lnot)P(t_{0},t_{1}), (\lnot)P_{\rho_{0}}(t)$.
In each case $ \dg(A)=1$.
\[
\infer*{\Gam_{end}:a_{1}}
{
\infer[(cut)]{\Gam,\Del;a\# b}
{
 \infer*{\Gam,\lnot A;a}{}
 &
 \infer*{A,\Del;b}{A,\Del_{0};1}
 }
}
\]
When $A$ is a $\Del_{0}$-formula, let $\mathcal{P}^{\prime}$ be the following with the false $\Del_{0}$-formula $\lnot A$ down to the end-sequent $\Gam_{end}$.
\[
\infer*{\Gam_{end},\lnot A: a_{1}^{\prime}}
{
\infer*{\Gam,\lnot A;a}{}
}
\]

Otherwise $A$ is a $P$-literal.
Eliminate the false literal $\lnot A$ by Lemma \ref{lem:delta0elim} to get the following $\mathcal{P}'$.
\[
\infer*{\Gam_{end}: a_{1}^{\prime}}
{
\infer*{\Gam;a}{}
}
\]
We claim that the proof $\calP'$ with the restricted stock assignment is a proof with stock, and
$a_{1}^{\prime}<a_{1}$.
Let $J$ be a rule $(D_{1})$ with $c={\sf c}(J)$ below the $(cut)$ with the cut formula $A$.
Then $d'=o(\Gam_{1}(,\lnot A);\calP',{\sf c})<o(\Gam_{1};\calP,{\sf c})=d$ for the upper sequent $\Gam_{1}$ of $J$.
Moreover we have $G_{D_{1}(c\oplus d)}(c\oplus \ome^{d'})\subset G_{D_{1}(c\oplus d)}(c\oplus \ome^{d})<c\leq c\oplus \ome^{d}$ 
by Proposition \ref{prp:oa}.\ref{prp:oa.1}, {\bf (p.2)} and {\bf (p2.1)}.
Hence by Proposition \ref{prp:G} we obtain
$D_{1}(c\oplus \ome^{d'})\in\calh_{c\oplus d}(D_{1}(c\oplus \ome^{d}))\cap\rho_{0}=D_{1}(c\oplus \ome^{d})$.

Next let $a_{1}=D_{0}(c_{0}\oplus a_{0})$ and $a_{1}^{\prime}=D_{0}(c_{0}\oplus a_{0}^{\prime})$,
where $a_{0}=o(\Gam_{0};\calP,{\sf c})$ and $a_{0}^{\prime}=o(\Gam_{0}(,\lnot A);\calP^{\prime},{\sf c})$
for the upper sequent $\Gam_{0}$ of the last rule $(D_{0})$.
Then we see $a_{0}^{\prime}<a_{0}$ from the above with $(D_{1})$, and from Proposition \ref{prp:oa}.\ref{prp:oa.2}
that
for each ordinal $\gam$ in the set $G_{a_{1}}(c_{0}\oplus a_{0}^{\prime})$,
either $\gam\in G_{a_{1}}(c_{0}\oplus a_{0})<c_{0}$ or
there exists a rule $(D_{1})_{\bet}\, J$ such that $\gam\in\{c\oplus \ome^{d'}\}\cup G_{a_{1}}(c\oplus \ome^{d'})$ 
with $c={\sf c}(J)$ and
$d'=o(\Gam_{1}(,\lnot A);\calP',{\sf c})$ for the upper sequent $\Gam_{1}$ of $J$.
We have $G_{a_{1}}(c\oplus \ome^{d'})\subset G_{a_{1}}(c\oplus \ome^{d})<c_{0}$ for $d=o(\Gam_{1};\calP,{\sf c})$.
On the other hand we have $c\oplus \ome^{d'}<c\oplus \ome^{d}\leq\bet_{0}<c_{0}$ for
$\bet=D_{0}(\bet_{0})$ by {\bf (p2)} and {\bf (p2.2)}.
Therefore $G_{a_{1}}(c_{0}\oplus a_{0}^{\prime})<c_{0}<c_{0}\oplus a_{0}$,
and we conclude $a_{1}^{\prime}<a_{1}$ from Proposition \ref{prp:G}.
\\

\noindent
{\bf Case 2}. The top is an axiom $(P_{\rho_{0}}\exi)$.

Let $C\equiv(\exi x[t<x\land P_{\rho_{0}}(x)])$.
Consider the uppermost and the lowest $(D_{1})$'s below the $(cut)$ whose cut formula is $C$.
We see that such a $(D_{1})$ exists below the cut from {\bf (h2)}.
\[
\infer[(D_{1})_{\bet}]{\cdots;D_{1}(c_{1}\oplus \ome^{b})}
{
 \infer*{\cdots}
 {
\infer[(D_{1})_{\alp}]{\Gam^{\prime}; b}
{
 \infer*{\Gam; b}
 {
  \infer[(cut)]{\Gam_{1},\Del_{1};b_{1}\# a_{1}}
  {
   \infer*{\Del_{1},\lnot C;b_{1}}{}
   &
   \infer*{C,\Gam_{1};a_{1}}
    {
    \infer[(P_{\rho_{0}}\exi)]{\Gam_{0},\exi x[t<x\land P_{\rho_{0}}(x)]; 1}{}
    }
  }
 }
}
}
}
\]
where 
there is no $(D_{1})$ above the $(cut)$ by {\bf (h6)}.
Let
$\ell=D_{1}(c_{1}\oplus 0)$.
We have $t\in\calh_{c_{1}}(D_{1}(c_{1}))\cap \rho_{0}=D_{1}(c_{1})=\ell$ by {\bf (p2.1)}
and Proposition \ref{prp:G}.
By inversions for the $A$-formula $\lnot C$, augmenting the sequent $\Gam_{1}$
 and eliminating false literals $t\not<\ell, \lnot P_{\rho_{0}}(\ell)$ we obtain
the following $\mathcal{P}^{\prime}$ with the new stock $c_{2}=c_{1}+1$ of the rules $(D_{1})$.

\[
\infer[(D_{1})_{\bet}]{\cdots;D_{1}( c_{2}\oplus \ome^{b^{\prime}})}
{
\infer*{\cdots}
{
 \infer[(D_{1})_{\alp}]{\Gam^{\prime}; b^{\prime}}
 {
  \infer*{\Gam; b^{\prime}}
  {
    \infer*[x:=\ell]{\Gam_{1},\Del_{1};b_{1}^{\prime}}{}
  }
 }
}
}
\]
Let us first check the condition {\bf (p2.1)} for the $(D_{1})_{\alp}$ in $\mathcal{P}^{\prime}$.
Any term occurring in $\calP'$ is in the closure of $\ell$ and terms occurring in $\calP$
under $+,\cdot,\ome$.
Hence it suffices to show $G_{D_{1}(c_{2}\oplus d)}(\ell)<c_{2}$,
which follows from $c_{1}<c_{2}$ and
$G_{D_{1}(c_{2}\oplus d)}(c_{1})<c_{1}$.

Next let us show $D_{1}(c_{2}\oplus \ome^{b^{\prime}})<D_{1}(c_{1}\oplus \ome^{b})$.
It is easy to see that $b_{1}^{\prime}\leq b_{1}$, and
$b^{\prime}+1<b$.
Moreover $G_{D_{1}(c_{1}\oplus \ome^{b})}(c_{2}\oplus \ome^{b^{\prime}})\subset G_{D_{1}(c_{1}\oplus \ome^{b})}(c_{1}\oplus \ome^{b})\cup G_{D_{1}(c_{1}\oplus \ome^{b})}(\ell)<c_{2}\leq c_{1}\oplus \ome^{b}$.
Hence by Proposition \ref{prp:G} we obtain $D_{1}(c_{2}\oplus \ome^{b^{\prime}})<D_{1}(c_{1}\oplus \ome^{b})$.

Finally we have $c_{2}\oplus b^{\prime}\leq c_{1}\oplus b<c_{0}$ for the stock $c_{0}$ of the last rule $(D_{0})$.
Hence {\bf (p2.2)} is enjoyed, and
$\calP'$ with the new stock is a proof with stock and $o(\calP')<o(\calP)$.
\\

\noindent
{\bf Case 3}. The top is an axiom $(P\exi)$.

First let $t\not<\ome_{1}$.
Since we are assuming that the end-sequent $\Gam_{end}$ is false,
the true literal $t\not<\ome_{1}$ vanishes at a $(cut)$.
$\calP$ be the following.
\[
\infer[(D_{0})_{\alp}]{\Gam_{end};D_{0}(c_{0}\oplus a)}
{
 \infer*{\Gam_{end};a}
 {
  \infer[(cut)]{\Gam_{1},\Del_{1};b_{0}\# a_{0}}
  {
   \infer*{\Del_{1},t<\ome_{1};b_{0}}{}
   &
   \infer*{t\not<\ome_{1},\Gam_{1};a_{0}}
    {
    \infer[(P\exi)]{\Gam_{0}, t\not<\ome_{1}, \exi x, y< \ome_{1}[t<x\land P(x,y)];1}{}
    }
  }
 }
}
\]
where the rule $(D_{0})_{\alp}$ is the last rule by {\bf (h7)}.
Eliminate the false $t<\ome_{1}$ to get the following $\calP'$.
\[
\infer[(D_{0})_{\alp}]{\Gam_{end};D_{0}(c_{0}\oplus a^{\prime})}
{
 \infer*{\Gam_{end};a^{\prime}}
 {
   \infer*{\Gam_{1},\Del_{1};b_{0}}{}
 }
}
\]
As in {\bf Case 1} we see that the resulting $\calP'$ is a proof with stock, and $D_{0}(c_{0}\oplus a^{\prime})<D_{0}(c_{0}\oplus a)$.

Next let $t$ be a closed term such that $t<\ome_{1}$, and
$C\equiv(\exi x, y< \ome_{1}[t<x\land P(x,y)])$.
\[
\infer[(D_{0})_{\alp}]{\Gam_{end};D_{0}(c_{0}\oplus a)}
{
 \infer*{\Gam_{end}; a}
 {
  \infer[(cut)]{\Gam_{1},\Del_{1};b_{1}\# a_{1}}
  {
   \infer*{\Del_{1},\lnot C;b_{1}}{}
   &
   \infer*{C,\Gam_{1};a_{1}}
    {
    \infer[(P\exi)]{\Gam_{0}(, t\not<\ome_{1}), \exi x, y< \ome_{1}[t<x\land P(x,y)]; 1}{}
    }
  }
 }
}
\]
where the rule $(D_{0})_{\alp}$ is the last rule by {\bf (h7)}.

Let $\alp\geq D_{0}(c_{0}\oplus a)>\ell=D_{0}(c_{0}\oplus 0)$.
Then $\ell>t$ 
by $t\in\calh_{c_{0}}(D_{0}(c_{0}\oplus 0))\cap D_{1}(0)=\ell$, {\bf (p2.1)}.
Let $s=F(c_{0}\oplus 0)$, i.e., $s=F_{\ell\cup\{\ome_{1}\}}(\rho_{0})$.

By inversions for the $A$-formula $\lnot C$ and eliminating false literals $\ell\not< \ome_{1}, s\not<\ome_{1}, t\not<\ell, \lnot P(\ell,s)$ we obtain the following with the new stock $c_{2}=c_{0}+1$ of the last rule $(D_{0})_{\alp}$.

\[
\infer[(D_{0})_{\alp}]{\Gam_{end};D_{0}(c_{2}\oplus a^{\prime})}
{
 \infer*{\Gam_{end}; a^{\prime}}
 {
   \infer*[x:=\ell, y:=s]{\Gam_{1},\Del_{1};b_{1}^{\prime}}{}
  }
 }
\]
The condition {\bf (p2.1)} for the $(D_{0})_{\alp}$ in $\mathcal{P}^{\prime}$ is seen to be fulfilled as in {\bf Case 2}.
{\bf (p2.1)} for rules $(D_{1})$ is enjoyed since
$G_{D_{1}(c\oplus d)}(\ell)=G_{D_{1}(c\oplus d)}(F(c_{0}\oplus 0))=\emptyset$ by $\ell,F(c_{0}\oplus 0)<\ome_{1}<D_{1}(c\oplus d)$
for any $c,d$.
As in {\bf Case 2} we see that $\calP'$ with the new stock is a proof with stock, and 
$D_{0}(c_{2}\oplus a^{\prime})<D_{0}(c_{0}\oplus a)$.

\subsection{top=rule}\label{subsec:toprule}

In this subsection we consider the cases when the top $\Phi$ is a lower sequent of 
one of explicit rules $(\lor),(\land),(\exi),(b\exi),(\fal), (b\fal)$ or $(Rfl),(ind)$ or one of implicit rules $(\lor),(\exi)$.
\\

\noindent 
{\bf Case 1}.
The top is the lower sequent of an explicit logical rule $J$.
Since the end-sequent consists solely of closed formulas, the main formula of $J$ is also closed.
\\

\noindent
{\bf Case 1.1}.
$J$ is a $(\fal)$:
Let $\mathcal{P}$ be the following.
\[
\infer[(D_{0})]{\Gam,\fal y\lnot A(y);D_{0}(c_{0}\oplus b)}
{
  \infer*{\Gam,\fal y\lnot A(y);b}
  {
    \infer[(\fal)\,J]{\Gam_{0},\fal y \lnot A(y); a_{0}+1}
   {
    \infer*{\Gam_{0},\fal y \lnot A(y),\lnot A(y); a_{0}}{}
  }
 }
}
\]
Note that the predicate $P_{\rho_{0}}$ does not occur in the end sequent $\Gam,\fal y\lnot A(y)$, and hence
any $(D_{1})$ does not change the descendants of the formula $\fal y \lnot A(y)$.
For the closed term $s\equiv\mu y.A(y)$ with $\fal y\lnot A(y)\lrarw \lnot A(s)$,
let $\mathcal{P}^{\prime}$ be the following.
\[
\infer[(D_{0})]{\Gam,\fal y\lnot A(y),\lnot A(s);D_{0}(c_{0}\oplus b^{\prime})}
{
\infer*{\Gam,\fal y\lnot A(y),\lnot A(s);b^{\prime}}
{
  \infer*[y:=s]{\Gam_{0},\fal y \lnot A(y),\lnot A(s);a_{0}^{\prime}}{}
 }
}
\]
where 
the closed term $s$ is substituted for the eigenvariable $y$.

Note that there is no rule $(D_{1})$ above the rule $(\fal)$
since no free variable occurs below $(D_{1})$ by {\bf (h1)}.
Let us check the condition {\bf (p2.1)} for a rule $(D_{i})\,(i=0,1)$ with its stock $c$ in $\mathcal{P}^{\prime}$.
Suppose that the formula $\fal y\lnot A(y)$ occurs in the upper sequent of $(D_{i})$.
Let $A(y)\equiv A(y;t_{1},\ldots,t_{k})$. Then
$G_{D_{i}(c\oplus d)}(f_{A}(t_{1},\ldots,t_{k}))\subset\bigcup_{m}G_{D_{i}(c\oplus d)}(t_{m})<c$ for $s\equiv\mu y.A(y)\equiv f_{A}(t_{1},\ldots,t_{k})$.
\\

\noindent
{\bf Case 1.2}.
$J$ is an $(\exi)$: 
\[
\infer*{\Gam_{end};a}
{
 \infer[(\exi)\,J]{\Gam_{0},\exi y\,A(y);a_{0}+1}
 {
  \Gam_{0},\exi y\,A(y),A(s);a_{0}
  }
 }
\]
where $s$ is a closed term.
If there is no rule affecting on descendant $\exi y\, A(y)$ of the main formula,
erase the rule $J$ to get a proof $\calP'$ of $\Gam_{end},A(s)$.
Suppose that there is a rule $J_{0}$ affecting on a descendant $\exi y\, A(y)$ of the main formula.
$J_{0}$ is one of the rules $(P\Sig_{1})$ and $(P_{\rho_{0}}\Sig_{1})$ 
since the predicate $P_{\rho_{0}}$ does not occur in the end-sequent $\Gam_{end}$.
\\

\noindent
{\bf Case 1.2.1}.
The rule $J_{0}$ is a $(P\Sig_{1})$:
Then $(\exi y\,A(y))\equiv(\vphi[\ome_{1},s_{0}])$ and
$(\exi y<t\,A^{\prime}(y))\equiv(\vphi^{t}[t_{0},s_{0}])$ for some closed terms $s_{0},t_{0}$.
\[
\infer*{\Gam_{end},\exi y<t\,A^{\prime}(y);a}
{
\infer[(P\Sig_{1})]{\Gam_{1},(\lnot P(t_{0},t),s_{0}\not<t_{0},) \vphi^{t}[t_{0},s_{0}];b}
{
 \infer*{\Gam_{1},\vphi[\ome_{1},s_{0}]; b}
 {
 \infer{\Gam_{0},\exi y\,A(y);a_{0}+1}
 {
  \Gam_{0},A(s);a_{0}
  }
 }
}
}
\]
If one of $\lnot P(t_{0},t)$ and $s_{0}\not<t_{0}$ is true, then eliminate 
one of the false literals $P(t_{0},t)$ and $s_{0}<t_{0}$
as in {\bf Case 3} of subsection \ref{subsec:topaxiom}.

Suppose that both $P(t_{0},t)$ and $s_{0}<t_{0}$ are true.
Then $\vphi[\ome_{1},s_{0}]\to\vphi^{t}[t_{0},s_{0}]$.
Let the closed false $\Del_{0}$-formula $A(s)$ go down to the end-sequent to get a proof
$\mathcal{P}^{\prime}$.
It is easy to see that $\calP'$ is a proof with stock such that $a^{\prime}<a$.
\[
\infer*{\Gam_{end},\exi y<t\,A^{\prime}(y),A(s); a^{\prime}}
{
\infer[(P\Sig_{1})]{\Gam_{1},(\lnot P(t_{0},t),s_{0}\not<t_{0},) \vphi^{t}[t_{0},s_{0}],A(s);b^{\prime}}
{
 \infer*{\Gam_{1},\vphi[\ome_{1},s_{0}],A(s); b^{\prime}}
 {
 \infer*{\Gam_{0},\exi y\,A(y),A(s); a_{0}}{}
 }
}
}
\]
{\bf Case 1.2.2}.
The rule $J_{0}$ is a $(P_{\rho_{0}}\Sig_{1})$.

Then $(\exi y\,A(y))\equiv(\vphi[s_{0}])$ and
$(\exi y<t\,A^{\prime}(y))\equiv(\exi y<t\,A(y))\equiv(\vphi^{t}[s_{0}])$ for a closed term $s_{0}$.
\[
\infer*{\Gam_{end},\exi y<t\,A(y); a}
{
\infer[(P_{\rho_{0}}\Sig_{1})]{\Gam_{1},(\lnot P_{\rho_{0}}(t),s_{0}\not<t,) \vphi^{t}[s_{0}]; b}
{
 \infer*{\Gam_{1},\vphi[s_{0}]; b}
 {
 \infer{\Gam_{0},\exi y\,A(y); a_{0}+1}
 {
  \infer*{\Gam_{0},\exi y\,A(y),A(s); a_{0}}{}
  }
 }
}
}
\]
If one of $\lnot P_{\rho_{0}}(t)$ and $s_{0}\not<t$ is true, then eliminate 
one of the false literals $P_{\rho_{0}}(t)$ and $s_{0}<t_{0}$
as in {\bf Case 3} of subsection \ref{subsec:topaxiom}.

Suppose that both $P_{\rho_{0}}(t)$ and $s_{0}<t$ are true.
Then $\vphi[s_{0}]\to\vphi^{t}[s_{0}]$.
Let the closed false $\Del_{0}$-formula $A(s)$ go down to the end-sequent to get a proof
$\mathcal{P}^{\prime}$.
It is easy to see that $\calP'$ is a proof with stock such that $a^{\prime}<a$.
\[
\infer*{\Gam_{end},\exi y<t\,A(y),A(s);a^{\prime}}
{
\infer[(P_{\rho_{0}}\Sig_{1})]{\Gam_{1},(\lnot P_{\rho_{0}}(t),s_{0}\not<t,) \vphi^{t}[s_{0}],A(s); b^{\prime}}
{
 \infer*{\Gam_{1},\vphi[s_{0}],A(s); b^{\prime}}
 {
 \infer{\Gam_{0},\exi y\,A(y),A(s); a_{0}}
 {
  \infer*{\Gam_{0},A(s); a_{0}}{}
  }
 }
}
}
\]
Other cases $(\lor),(\land), (b\fal),(b\exi)$ are similar.
\\

\noindent
{\bf Case 2}.
The top is the lower sequent of a $(Rfl)$:
Let $A(x)\equiv(\exi z\exi w[z\in P_{\rho_{0}}\land B(x)])\, (B\in\Del_{0})$,
$A^{(y)}(x)\equiv(\exi z<y\exi w< y[z\in P_{\rho_{0}}\land B])$.
\[
\infer[J]{\Del; D_{1}(c_{1}\oplus \ome^{a})}
{
 \infer*{\Del_{1}}
 {
  \infer[J_{1}]{\Del_{2}^{\prime}; a}
  {
   \infer*{\Del_{2}; a}
   {
    \infer[(Rfl)]{\Gam; a_{0}\# a_{1}}
     {
     \Gam, \fal x< t\, A(x); a_{0}
     &
     t\not< y, \exi x< t \lnot A^{(y)}(x),\Gam; a_{1}
     }
    }
  }
 }
}
\]
where 
$J_{1}$ is the uppermost $(D_{1})_{\alp_{1}}$ and $J$ is the lowermost $(D_{1})_{\alp}$ below the $(Rfl)$.
Such a $(D_{1})$ exists by {\bf (h5)}.
Let $\mathcal{P}'$ be the following.
\[
\infer[(cut)]{\Del; \ell+r}
{
 \infer{\Del,\fal x< t A^{(\ell)}(x); D_{1}(c_{1}\oplus \ome^{a_{\ell}})}
 {
  \infer*{\Del_{1},\fal x< t A^{(\ell)}(x); a_{\ell}}
  {
   \infer{\Del_{2}^{\prime},\fal x< t A^{(\ell)}(x); a_{\ell}}
   {
    \infer[(D_{1})_{\ell}]{\Del_{2},\fal x< t A^{(\ell)}(x); a_{\ell}}
    {
     \infer*{\Del_{2},\fal x< t A(x);a_{\ell}}
      {\Gam,\fal x< t A(x); a_{0}}
      }
     }
   }
  }
&
 \infer[J_{r}]{\exi x< t\lnot A^{(\ell)}(x),\Del; D_{1}(c_{2}\oplus \ome^{a_{r}})}
 {
  \infer*{\exi x< t\lnot A^{(\ell)}(x),\Del_{1}; a_{r}}
  {
   \infer[(D_{1})_{\alp_{1}}]{\exi x< t\lnot A^{(\ell)}(x),\Del_{2}^{\prime}; a_{r}}
   {
    \infer*{\exi x< t\lnot A^{(\ell)}(x),\Del_{2};a_{r}}
      {
       \infer*[y:=\ell]{\exi x< t\lnot A^{(\ell)}(x),\Gam; a_{1}^{\prime}}{}
     }
    }
  }
 }
}
\]
where $\ell:=D_{1}(c_{1}\oplus \ome^{a_{\ell}})$, $r:=D_{1}(c_{2}\oplus \ome^{a_{r}})$,
 and the stock of the rule $J_{r}$ is enlarged to $c_{2}=c_{1}\#\ome^{a_{\ell}}+1$.
We see $t<\ell$ from {\bf (p2.1)}, and the false literal $t\not<\ell$ is eliminated.

In $\mathcal{P}$,  $h(\Del)\geq \dg(\exi x< t \lnot A^{(y)}(x))= \dg(\fal x< t A^{(\ell)}(x))$ by {\bf (h5)}.
Thus the introduced $(cut)$ in $\mathcal{P}'$ enjoys {\bf (h3)}.
There is no $(D_{1})$ above the $(Rfl)$ by {\bf (h6)}.
In the left part of the $(cut)$, a new $(D_{1})_{\ell}$ arises with its stock $c_{1}$,
cf.\,{\bf (p2)}.
In the upper sequent of the right rule $(D_{1})_{\alp_{1}}$, 
a bounded sentence $\exi x< t\lnot A^{(\ell)}(x)$ is added, cf.\,the definition of the rule $(D_{1})$.
For the condition {\bf (p2.1)} of the right rule $(D_{1})_{\alp_{1}}$ we have
$G_{D_{1}(c_{2}\oplus d)}(\ell)\subset \{c_{1}\oplus \ome^{a_{\ell}}\}\cup G_{D_{1}(c_{2}\oplus d)}(c_{1}\oplus \ome^{a_{\ell}})\subset \{c_{1}\oplus \ome^{a_{\ell}}\}\cup G_{D_{1}(c_{2}\oplus d)}(c_{1}\oplus \ome^{a})<c_{2}$.
It is clear that $\ell<D_{1}(c_{1}\oplus \ome^{a})$.
On the other hand we have $c_{2}\oplus \ome^{a_{r}}= c_{1}\#\ome^{a_{\ell}}\#\ome^{a_{r}}+1<c_{1}\oplus \ome^{a}$, and
$G_{r}(c_{2}\oplus \ome^{a_{r}})\subset G_{r}(c_{1}\oplus \ome^{a})\cup G_{r}(\ell)<c_{2}<c_{1}\oplus \ome^{a}$.
Thus $r<D_{1}(c_{1}\oplus \ome^{a})$.
From this we see that $\calP'$ is a proof with stock, and $o(\calP')<o(\calP)$.
\\

\noindent
{\bf Case 3}.
The top is the lower sequent of an $(ind)$.
\[
\infer[(D_{0})]{\cdots; D_{0}(c_{0}\oplus b_{0})}
{
\infer*{\cdots;b_{0}}
{
\infer[(D_{1})]{\cdots;c_{1}\oplus b_{1}}
{
\infer*{\cdots;b_{1}}
{
\infer[(ind)]{(s\not< t,) \Gam;a}{\Gam,\lnot\fal x< y A(x), A(y);a_{0} & \Gam,\lnot A(s);a_{1}}
}
}
}
}
\]
where $(D_{1})$ is the uppermost one.
Such a $(D_{1})$ exists by {\bf (h4)}.
There is no $(D_{1})$ nor $(ind)$ above the $(ind)$ by {\bf (h6)} and {\bf (h4)}.
By {\bf (p1)} we have $ \dg(A(y))=a_{1}$, $a_{0}<\ome$ and
$a=(a_{0}+a_{1}+2)\times t$ for the closed term $t$.
\\
{\bf Case 3.1}. $s\not< t$:
Then the true literal $s\not< t$ remains in the lower sequent.
Eliminate the false literal $s<t$.
\\

\noindent
{\bf Case 3.2}. $s< t$:
Assuming $\lnot A(s)$ is an $\exi$-formula, 
let $P'$ be the following: 
{\footnotesize
\[\hspace*{-7mm}
\infer[(D_{0})]{\cdots; D_{0}(c_{0}\oplus b_{0}')}
{
\infer*{\cdots; b_{0}'}
{
\infer[(D_{1})]{\cdots; c_{1}\oplus b_{1}'}
{ 
\infer*{\cdots;b_{1}'}
{
 \infer[(cut)]{(s\not< t,) \Gam;a'}
 {
  \infer[(b\fal)]{\Gam,\fal x< s A(x)}
   {
    \infer[(ind)]{\Gam,y\not< s,A(y);(a_{0}+a_{1}+2)\times s}
     {
      \infer*{\Gam,\lnot\fal x< y A(x),A(y);a_{0}}{}
      &
       \infer*[P(A)]{\Gam,\lnot A(y),A(y);a_{1}}{}
      }
    }
   &
    \infer[(cut)]{\lnot\fal x< s\, A(x),\Gam}
     {
      \infer*[y:=s]{\Gam,\lnot\fal x< s\, A(x),A(s);a_{0}}{}
      &
      \infer*{\lnot A(s),\Gam;a_{1}}{}
     }
  }
 }
}
}
}
\hspace*{-10mm}
P'
\]
}
where $P(A)$ denotes a proof of $\Gam,\lnot A(y),A(y)$ which is canonically constructed from logical inferences, 
cf.\,Tautology lemma \ref{lem:depth}.

We have $h_{0}(\Gam)\geq \dg(\fal x< a A(x))\geq \dg(A(a))$ by {\bf (h4)}, and hence {\bf (h3)} holds for the introduced $(cut)$'s.
Also $a'=(a_{0}+a_{1}+2)\times s+a_{0}+a_{1}+1<(a_{0}+a_{1}+2)\times t=a$.
Since no essentially new term is created here, {\bf (p2.1)} is fulfilled with $\mathcal{P}^{\prime}$.

If $\lnot A(s)$ is not an $\exi$-formula, then upper sequents of the upper cut should be interchanged.
Note that $a_{1}+a_{0}=a_{0}+a_{1}$ for $a_{0},a_{1}<\ome$:
\[
\infer{(s\not<t, ) \Gam}
{
 \infer*{\Gam,\fal x< s A(x)}{}
&
 \infer{\lnot \fal x< s A(x),\Gam}
 {
  \infer*{\lnot A(s),\Gam;a_{1}}{}
 &
  \infer*{\Gam,\lnot\fal x< s A(x),A(s);a_{0}}{}
 }
}
\msten P'
\]
{\bf Case 4}.
The top $\Phi$ is the lower sequent of a logical rules $(\lor), (\exi), (b\exi)$.
Consider the cases when the logical rule is one of $ (\exi),(b\exi)$, which is denoted $(\exi)$.
The case $(\lor)$ is similar.
Let the main formula of the logical rule be a formula $\exi x< t \, A(x)$ with a minor formula $A(s)$,
where $t$ denotes either a term or $\rho_{0}$, $(\exi x<\rho_{0}\,A(x)):\equiv(\exi x\, A(x))$.
Let $J$ denote the $(cut)$ at which the descendant $\exi x< t'\, A'(x)$ of $\exi x<t\, A$ vanishes.
\\

\noindent
{\bf Case 4.1}.
$\exi x< t'\, A'(x)$ is a $\Del_{0}$-formula:
Let $\mathcal{P}$ be the following.
\[
 \infer*{\Gam_{end};c}
 {
  \infer[(cut)\, J]{\Gam,\Del;a\# b}
  {
   \infer*{\Gam,\lnot \exi x< t^{\prime} \, A^{\prime}(x);a}{}
   &
   \infer*{\exi x< t^{\prime} \, A^{\prime}(x),\Del;b}
   {
    \infer[(\exi)]{\exi x< t\, A(x),\Del_{0};b_{0}+1}{\exi x< t\, A(x),A(s),\Del_{0};b_{0}}
   }
  }
 }
\]

One of $\lnot \exi x<t^{\prime}\, A^{\prime}(x),\exi x<t^{\prime}\, A^{\prime}(x)$ is false.
When $\exi x<t^{\prime}\, A^{\prime}(x)$ is false, let the false $\Del_{0}$-formula $\exi x<t^{\prime}\, A^{\prime}(x)$
go down to the end-sequent.
\[
 \infer*{\Gam_{end},\exi x<t^{\prime}\, A^{\prime}(x);c(b)}
 {
   \infer*{\exi x<t^{\prime}\, A^{\prime}(x),\Del;b}
      {
    \infer[(\exi)]{\exi x< t\, A(x),\Del_{0};b_{0}+1}{\exi x< t\, A(x),A(s),\Del_{0};b_{0}}
   }
 }
\]

When $\exi x<t^{\prime}\, A^{\prime}(x)$ is true, we are in {\bf Case 1.1} of this subsection.
\[
 \infer*{\Gam_{end},\lnot \exi x<t^{\prime}\, A^{\prime}(x);c(a)}
 {
  \infer*{\Gam,\lnot\exi x<t^{\prime}\, A^{\prime}(x);a}{}
  }
\]
In each case we have $c(b),c(a)<c$.
In what follows assume that $\exi x< t'\, A'(x)$ is not a $\Del_{0}$-formula
\\

\noindent
{\bf Case 4.2}.
The descendant $\exi x<t^{\prime}\,A^{\prime}(x)$ 
may differ from the main formula $\exi x<t\,A(x)$ due to a rule $(D_{1})$
with $t=\rho_{0}$
when either $(\exi x<t\,A(x))\equiv(\exi x\exi w[x\in P_{\rho_{0}}\land B(x,w)])\, (B\in\Del_{0})$,
or
$(\exi x<t\,A(x))\equiv(\exi x[s\in P_{\rho_{0}}\land B(s,x)])$.
The case when a rule $(P\Sig_{1}), (P_{\rho_{0}}\Sig_{1})$ change a descendant of the main formula is excluded
since we are assuming that $\exi x< t'\, A'(x)$ is not a $\Del_{0}$-formula.
Note that there is no $(D_{1})$ nor $(D_{0})$ above the $(cut)\, J$ by {\bf (h6)}, and there is a $(D_{0})$ below the vanishing cut by {\bf (h7)}.
Since $\exi x<t\,A'(x)$ is not a $\Del_{0}$-formula,
$ \dg(\exi x< t \, A'(x))>0$,
and there exists an $(h)$ below the vanishing cut by {\bf (h3)}.
Consider the uppermost $(h)$.
\[
\infer[(h)]{\Lam;\ome^{c}}
{
 \infer*{\cdots;c}
 {
  \infer[(cut)]{\Gam,\Del;a\# b}
  {
   \infer*{\Gam,\lnot \exi x< t \, A'(x);a}{}
   &
   \infer*{\exi x< t \, A'(x),\Del;b}
   {
    \infer[(\exi)]{\exi x< t\, A(x),\Del_{0};b_{0}+1}{\exi x< t\, A(x),A(s),\Del_{0};b_{0}}
   }
  }
 }
}
\]
Since $h_{0}(\Gam,\Del)\geq \dg(\exi x< t \, A'(x))> \dg(A'(s))$, we have 
$h_{0}(\Lam)=h_{0}(\Gam,\Del)-1\geq \dg(A'(s))$ for {\bf (h3)}.
Assuming that $\lnot A(s)$ is an $E$-formula, let $\calP'$ be the following.

\[
\infer[(cut)]{\Lam; \ome^{c_{\ell}}\#\ome^{c_{r}}}
{
 \infer[(h)]{\Lam,A(s); \ome^{c_{\ell}}}
 {
  \infer*{\cdots;c_{\ell}}
  {
   \infer{\Gam,\Del,A(s);a\# b'}
   {
    \infer*{\Gam,\lnot \exi x< t \, A'(x);a}{}
    &
    \infer*{\exi x< t\, A'(x),A(s)\Del;b'}
    {
     \exi x< t\, A(x),A(s),\Del_{0};b_{0}
    }
   }
  }
 }
&
\infer[(h)]{\lnot A(s),\Lam;\ome^{c_{r}}}
{
 \infer*{\cdots;c_{r}}
 {
   \infer*[x:=s]{\lnot A'(s),\Gam;a^{\prime}}{}
  }
 }
}
\]
Note that there may occur a $(D_{1})$ above the left part of the $(cut)$ in $\mathcal{P}$.
Let $(D_{1})_{\alp}$ be a rule occurring above the left upper sequent of the $(cut)$ such that
its lower sequent contains an ancestor $\lnot\exi x<t\,A'(x)$ of the left cut formula.
We have to verify the condition {\bf (p2.1)} for the $(D_{1})$ in $\mathcal{P}^{\prime}$.
Let $c_{1}$ be the stock of the $(D_{1})_{\alp}$.
Then $G_{D_{1}(c_{1}\oplus d)}(t)<c_{1}$ for any $d$, where
$t<\rho_{0}$ since 
a formula $\lnot\exi x<t\,A^{\prime}(x)$ is in the upper sequent of the $(D_{1})_{\alp}$,
either $A^{\prime}\equiv A$ or $(A^{\prime})^{(\alp)}\equiv A$,
and there occurs no implicit formula with an unbounded universal quantifier in an upper sequent of a rule $(D_{1})$
by the definition of the rule.
Hence $s<t\in\calh_{D_{1}(c_{1}\oplus d)}(c_{1}\oplus d)\cap\rho_{0}=D_{1}(c_{1}\oplus d)$,
and $G_{D_{1}(c_{1}\oplus d)}(s)=\emptyset$.This shows {\bf (p2.1)}.

This completes a proof of Lemma \ref{lem:reflect}, and Theorem \ref{th:main}.\ref{th:main1}.
Theorem \ref{th:main}.\ref{th:main2} is seen from the proof of Lemma \ref{lem:reflect}
by restricting to a subset $T(\ome_{k}(\rho_{0}+1))$ and Lemma \ref{lem:lowerbndreg}.


\begin{thebibliography}{99}





\bibitem{ptMahlo} T. Arai, Proof theory for theories of ordinals I:recursively Mahlo ordinals, 
Ann. Pure Appl. Logic 122 (2003), 1-85.

\bibitem{liftupK}T. Arai, Proof theory of weak compactness, Jour. Math. Logic 13(2013), 1350003, 26pages

\bibitem{liftupZF}T. Arai, 
Lifting proof theory to the countable ordinals: Zermelo-Fraenkel's set theory,
Jour. Symb. Logic 79 (2014), 325-354.


\bibitem{IntSet}T. Arai, Intuitionistic fixed point theories over set theories, 
Arch. Math. Logic 54 (2015), 531-553.

\bibitem{Buchholz86} W. Buchholz,
A new system of proof-theoretic ordinal functions,
Ann. Pure Appl. Logic 32 (1986), 195-208.

\bibitem{Buchholz92} W. Buchholz, A simplified version of local predicativity, P. H. G. Aczel, H. Simmons and S. S. Wainer(eds.), Proof Theory, Cambridge UP, 1992, pp. 115-147.

\bibitem{Gentzen38}
G. Gentzen,
Neue Fassung des Widerspruchsfreiheitbeweis f\"ur die reine Zahlentheorie,
Forschungen Zur Logik und zur Grundlegung der exakten Wissenschaften, Neue Folge 4(1938), 19-44.


\bibitem{Takeuti67}
G. Takeuti,
Consistency proofs of subsystems of classical analysis,
Ann.  Math., 86(1967), 299-348.



\end{thebibliography}
\end{document}